\documentclass[reqno,11pt]{amsart}
\usepackage{amsthm,amsmath,amsfonts,amssymb,amscd,mathrsfs,graphics,diagrams}
\usepackage{txfonts}%,pxfonts}
\usepackage{hyperref,supertabular}

\DeclareMathOperator{\LG}{{LG}}
\DeclareMathOperator{\CY}{{CY}}
\DeclareMathOperator{\PD}{{PD}}
\renewcommand{\top}{{top}}
\newcommand{\C}{{\mathbb{C}}}

\newcommand{\D}{{\mathscr D}}

\newcommand{\E}{{\mathscr E}}
\newcommand{\Eo}{{\E_{\ominus}}}
\newcommand{\End}{{\operatorname{End}}}

\newcommand{\I}{{\mathcal I}}
\newcommand{\Io}{{\I_{\ominus}}}

\newcommand{\K}{{\mathcal K}}

\newcommand{\N}{{\mathbb N}}

\newcommand{\R}{{\mathbb{R}}}

\newcommand{\Res}{{\operatorname{Res}}} %residue

\newcommand{\Z}{{\mathbb{Z}}}

\newcommand{\ch}{{\mathscr H}}
\newcommand{\cho}{{\ch_{\ominus}}}

\newcommand{\hch}{{\widehat{\mathscr H}}}

\newcommand{\ib}{{\bar{i}}}

\newcommand{\jb}{{\bar{j}}}
\newcommand{\kb}{{\bar{k}}}

\newcommand{\mb}{{\bar{m}}}
\newcommand{\mub}{{\bar{\mu}}}

\newcommand{\res}{\operatorname{res}}
\newcommand{\pat}{{\partial}}
\newcommand{\bpat}{{\bar{\partial}}}

\newcommand{\Rb}{\underline{r}}
\newcommand{\Rbb}{\underline{R}}

\newtheorem{thm}{Theorem}[section]
\newtheorem{lm}[thm]{Lemma}
\newtheorem{prop}[thm]{Proposition}
\newtheorem{crl}[thm]{Corollary}

\theoremstyle{definition}
\newtheorem{rem}[thm]{Remark}

\newtheorem{df}[thm]{Definition}

\theoremstyle{remark}
 
\newtheorem{defthm}[thm]{{\bf Definition-Theorem}}

\title{Constructing the LG/CY isomorphism between $tt^*$ geometries}
\date{}
\author{Huijun Fan$^\dag$}
\address{School of Mathematical Sciences, Peking University, Beijing, China}
\email{fanhj@math.pku.edu.cn}
\author{Tian Lan}
\address{School of Mathematical Sciences, Peking University, Beijing, China}
\email{lantian8741@gmail.com}
\author{Zongrui Yang$^\dag$}
\address{Department of Mathematics, Columbia University, New York, NY, USA}
\email{zy2417@columbia.edu}
\thanks{$^\dag$ The first author is supported by NSFC(11271028, 11325101, 11671033, 11831017,11890660, 11890661) and NSFC-RFC 1201101428, NSFC-RFBR 11661131005. $^\dag$The third author is partially supported by Ivan Corwin's NSF grant DMS: 1811143 as well as the Fernholz Foundation's "Summer Minerva Fellows" program.}

\begin{document}
\maketitle
%\tableofcontents

\begin{abstract} For a  nondegenerate homogeneous polynomial  $f\in\mathbb{C}[z_0, \dots, z_{n+1}]$ with degree $n+2$, we can obtain a $tt^*$ structure from the Landau-Ginzburg model $(\C^{n+2}, f)$ and a (new) $tt^*$ structure on the Calabi-Yau hypersurface defined by the zero locus of $f$ in $\C P^{n+1}$. We can prove that the big residue map considered by Steenbrink gives an isomorphism between the two $tt^*$ structures. We also build the correspondence for non-Calabi-Yau cases, and it turns out that only partial structure can be preserved. As an application, we show that the $tt^*$ geometry structure  of Landau-Ginzburg model on relavant deformation space uniquely determines the $tt^*$ geometry structure on Calabi-Yau side. This explains the folklore conclusion in physical literature. This result is based on our early work \cite{FLY}. 
\end{abstract}

\section{Introduction}

Mirror symmetry was found by physicists when they studied the string theory in 1980's. B. Greene and R. Plesser \cite{GP} observed firstly a strange duality of Hodge numbers between two T-dual Calabi-Yau 3-folds $M$ and $\check{M}$: $h^{3-p,q}(M)=h^{p,q}(\check{M})$ and $h^{1,1}(M)=h^{2,1}(\check{M})$. Furtheremore, by the work of P. Candelas, C. Xinia, P. S. Green and L. Parkes (ref. \cite{Y}), people have realized a mysterious duality between the two topological field theories defined respectively on $M$ and $\check{M}$. Note that $h^{1,1}(M)$ is the dimension of the deformation space of K\"{a}hler (symplectic) structures on $M$, while $h^{2,1}(\check{M})$ corresponds to the dimension of the deformation space of complex structures on $\check{M}$. Usually we call the geometrical structure related to the symplectic structure as A-model, and the corresponding topological field theory is called the A-theory. On the other hand, we call the geometrical structure related to the complex structure as B-model, and the corresponding topological field theory is called the B-theory. The two topological field theories come from the A or B twists of the quantum theory of the Calabi-Yau nonlinear sigmal models.  

The A theory of a Calabi-Yau manifold can be formulated as the Gromov-Witten theory in mathematics. However, the widely accepted integral B-theory of a Calabi-Yau manifold for higher genus has not been rigorously built, despite of many attempts (ref. \cite{BCOV, HKQ, CL}). The genus $0$ part of the B-theory was completely understood as a \textquotedblleft special geometry \textquotedblright(\cite{St, Du}) based on the earlier work of Cecotti and Cecotti-Vafa (\cite{Ce1, Ce2, CV1}) . This special geometry was studied by C. Hertling \cite{Het1} and was called $tt^*$-geometry. The mirror symmetry conjecture between Calabi-Yau manifolds for genus $0$ was solved in many cases (ref.\cite{HKKPTVVZ}), and the proof for genus $1$ and $2$ cases for quintic Calabi-Yau cases can be found in \cite{Zi,GJR}.  

In addition to the Calabi-Yau nonlinear sigma model, there is another physical model in supstring theory, called the Landau-Ginzburg (LG) model. The geometry of the LG model consists of a noncompact K\"{a}hler manifold with a holomorphic function (called superpotential) defined on it. For a supersymmetric LG model, there are also A and B type topological field theories. The A-theory of a LG model has been constructed by T. Jarvis, Y. Ruan and the first author based on Witten's r-spin theory (ref. \cite{FJR}) and is called FJRW theory now. The B-side of a LG model is closely related to the singularity theory and can be partially described by Saito's Frobenius manifold structure \cite{S1, S3, S4} and Givental's quantization method in the case that the symmetry group of the superpotential function is trivial. Hence combining the A and B models of the Calabi-Yau and LG models, there is a global mirror symmetry picture (ref. \cite{CIR}) and can be briefly described by the following diagram:
\begin{equation}\label{intr-mirr-diag}
\begin{diagram}
\text{LG A theory} &\rTo^{mirror} & \text{LG B theory}\\
\dTo^{\text{\small LG/CY corresp.}}& &\dTo_{\text{\small LG/CY corresp.}}\\
\text{CY A theory} &\rTo^{mirror} & \text{CY B theory}
\end{diagram}
\end{equation}

Note that the vertical line is given by the LG/CY correspondence conjecture proposed by physicists (ref.\cite{FJR}). The A-model conjecture has been studied extensively (see \cite{CR1, CR2, GS, FJR2, CFGKS} and references there). 

Unlike the A-model theory, the B-theory should contain not only the holomorphic part but also the anti-holomorphic part, which forms the so called $tt^*$-geometrical structures. Since Saito's Frobenius manifold structures or the categoric LG/CY correspondence proved by D. Orlov \cite{Or} only concerns the holomorphic part, one needs to add the real structure into the pictures. This was the motivation to study the $tt^*$-geometrical structures. The $tt^*$-geometrical structures for A-model have been considered by H. Iritani \cite{I}. In this paper, we only concern the $tt^*$-geometrical structures for B-models. The $tt^*$-geometries for Calabi-Yau manifolds appeared in the papers \cite{Du, BCOV, St}. In 2011, the first author \cite{Fan} provided another approach to build the $tt^*$ geometrical structures for LG models by considering the variation of Hodge structures related to the twisted Cauchy-Riemann operators. This is much like the way to get the VHS for compact K\"ahler manifolds. This approach and its various applications of $tt^*$ geometrical structures have been furtherly developed in the work \cite{Wen,T, LW}. 

The $tt^*$ geometry has also been extensively studied in many physical literatures (ref. \cite{CV2, BC, AB, AB2, AB3} and etc.). These papers took 3-dimensional CY hypersurfaces as examples and considered the computations of period integrals and Weil-Pertersson metrics on the moduli spaces via the  corresponding LG model. We notice that the consideration of the LG/CY correspondence between Hodge structures and many ideas appeared very early (ref. mathematical literatures \cite{CG} and physical literature \cite{Ce1, Ce2, CV1}). However, in this paper we build this correspondence for any dimensional Calabi-Yau hypersurface in projective space based in mathematical rigor, which completely solved this problem appeared in our early work \cite{FLY}. The main proof of this paper appeared in the PhD thesis of the second author \cite{L}, which is based on the careful study of the relations between various residue maps, in particular Steenbrink's big residue maps \cite{Ste}.   

To formulate the main result of this paper, let us introduce the definition of $tt^*$ geometry (ref. \cite{Fan, Wen,T}), which is defined as a category consisting of objects as $tt^*$ bundles (or $tt^*$ structures), and morphisms as the embeddings. 

\begin{df}[$tt^*$ bundle]
A $tt^*$ bundle $\E=(H\rightarrow M,\kappa,\eta,D,\bar{D}, C,\bar{C})$ consists of the following data: 
\begin{itemize}
\item $H\rightarrow M$ is a complex vector bundle (called the Hodge bundle);

\item $\kappa: H\rightarrow H$ is a complex anti-linear involution, i.e. $\kappa^2=Id, \kappa(\lambda \alpha)=\bar{\lambda}\kappa(\alpha), \forall\lambda\in \mathbb{C}$ , $\forall \alpha\in\Gamma(H)$ ($\kappa$ is called a real form);

\item $\eta$ is a nondegenerate pairing on $H$ and together with the real form $\kappa$ induces a Hermitian metric $g(u,v)=\eta(u, \kappa v)$, $\forall u,v\in\Gamma(H)$ (called $tt^*$ metric);

\item $\D=D+\bar{D}+C+\bar{C}$ is a flat connection on $H$ such that $D+\bar{D}$ is the Chern connection of $g$ (w. r. t. the holomorphic structure given by $\bar{D}$) and $C$ and $\bar{C}$ are $C^{\infty}$(M)-linear maps
\begin{displaymath}
C:C^{\infty}(H)\rightarrow C^{\infty}(H)\otimes \Lambda^{1,0}(M),\; \bar{C}:C^{\infty}(H)\rightarrow C^{\infty}(H)\otimes \Lambda^{0,1}(M)
\end{displaymath}
satisfying
\begin{enumerate}
\item $g$ is real with respect to $\kappa$, i.e. $g(\kappa(u),\kappa(v))=\overline{g(u,v)}$.
\item $(D+\bar{D})\kappa=0, \bar{C}=\kappa\circ C\circ\kappa$.
\item $\bar{C}$ is the adjoint of $C$ with respect to $g$, i.e. $g(C_Xu,v)=g(u,\bar{C}_{\bar {X}}v)$ $\forall u,v\in\Gamma(H), \forall X\in\Gamma(T^{1,0}(M))$.
\end{enumerate}
The operators $C, \bar{C}$ are called the Higgs fields and the connection $\D$ is called the Gauss-Manin connection.We denote $\nabla=\D^{1,0}=D+C$ and $\bar{\nabla}=\D^{0,1}=\bar{D}+\bar{C}$.  
\end{itemize}
\end{df}
Notice that the vector bundle $H\rightarrow M$ can be equipped with more than one holomorphic structures, for example, both $\nabla^{(0,1)}$ and $D^{(0,1)}$ can define a holomorphic structure on $H\rightarrow M$.

Then we give the definition of the embeddings of $tt^*$ bundles.
\begin{df} Let $\E_i=(H_i\longrightarrow M_i,\kappa_i,\eta_i,D_i,C_i,\bar{C_i}), i=1,2,$ be two $tt^*$ bundles. An embedding $\Phi=(\phi, \phi^{\prime})$ of  two holomorphic bundles
\begin{diagram}
H_1&\rTo^{\phi^{\prime}}& H_2\\
\dTo& & \dTo\\
M_1& \rTo^{\phi}& M_2
\end{diagram}
is called an embedding from the $tt^*$ bundle $\E_1$ to $\E_2$ if $\phi$ and $\phi^{\prime}$ are holomorphic maps and the following hold: $\forall p\in M_1,  X\in T_p^{1,0} M_1, u, v\in (H_1)_p$,
\begin{enumerate}
\item $\eta_1(u, v)=\eta_2\circ \phi (\phi^{\prime}(u), \phi^{\prime}(v))$.
\item $\kappa_2\circ \phi^{\prime}=\phi^{\prime}\circ \kappa_1$.
\item $\phi^{\prime}((D_1)_X u)=(D_2)_{\phi_*(X)}(\phi^{\prime}(u))$ and $\phi^{\prime}\circ \kappa_1((\bar{D}_1)_{\bar{X}} (\kappa_1(u)))=(\bar{D}_2)_{\overline{\phi_*({X})}}(\kappa_2(\phi^{\prime}(u))$
\item $\phi^{\prime}((C_1)_X u)=(C_2)_{\phi_*(X)}(\phi^{\prime}(u))$
\end{enumerate}
Moreover, if $\Phi$ is a bundle isomorphism, we say that $\Phi$ is an isomorphism between $tt^*$ bundles.
\end{df}

After introducing the definition of  $tt^*$ bundles, we describe the two types of $tt^*$ geometries built for the Calabi-Yau model and the LG model in the consequent chapters. 

In Section 2, we will review the construction of the (big)  $tt^*$ bundles for LG models in \cite{Fan}, the small  $tt^*$ bundles introduced in \cite{FLY} and the (old) $tt^*$ bundles for the Calabi-Yau models. Theorem \ref{fro-iso} gives the correspondence between the small  $tt^*$ bundles for LG model and the (old)  $tt^*$ bundles for the Calabi-Yau models except the real structures. This result was proved in \cite{FLY} by considering the small residue map $r'$ (ref.(\ref{sec2:defn-smal-resi-1})) introduced by Carlson and Griffiths in \cite{CG}.

Section 3 concerns the relations between various topological residue maps. A big residue map $\Rbb$ acting on the whole Milnor ring of the superpotential function was introduced by Steenbrink in \cite{Ste} by considering the compactification of Milnor fibers in a projective space. It turns out that the small residue map can be factorized through the big residue map and only the elements with appropriate degree can be nonzero after acting by the small residue map. This explains why the small residue map is defined as only acting on the subring of the Milnor ring.  

Section 4 builds the complete LG/CY correspondence between the small $tt^*$ bundles for a LG model and the new $tt^*$ bundle for the corresponding Calabi-Yau model. Note that in the previous result, Theorem \ref{fro-iso}, the small residue map does not preserve the real structures. This forces us to change the known structures. Our method is to replace the small residue map by the big residue map, which can be defined topologically. This implies that the real structures commute with the big residue map. Besides that the Gauss-Manin connections can also commute with the big residue maps. Therefore, we can pushforward the Higgs fields in the LG model to the Calabi-Yau model by the big residue map to get the new Higgs fields which are different to the (old) Higgs fields coming from the Griffiths' transversality theorem. This introduces a new  $tt^*$ bundle for the Calabi-Yau model. Finally, we can normalize the pairings in the two sides they have the exact pullback relations. The above conclusions will be proved in Section 4 and we have the main theorem of this paper.   

\begin{thm}[Theorem \ref{complete-corr}]\label{main-1}
Let $f:\C^{n+2}\rightarrow \C$ be a nondegenerate homogeneous polynomial of degree $(n+2)$. Let 
$$
\E^{\LG}=({H}^{\LG}\to M, \kappa^{\LG},\eta^{\LG}, D^{\LG}, \bar{D}^{\LG}, C^{\LG}, \bar{C}^{\LG})
$$
be the small  $tt^*$ bundle for the LG model given in Definition-Theorem \ref{thm-small-lg} and let 
$$
\E^{\CY}=(H^{\CY}\to M, \kappa^{\CY},\eta^{\CY}, D^{\CY}, \bar{D}^{\CY}, C^{\CY}, \bar{C}^{\CY})
$$
be the new  $tt^*$ bundle, then the map $\Rbb \circ \Phi: \E^{\LG}\to \E^{\CY}$ is an isomorphism, where the map $\Phi$ is defined in Theorem \ref{sec2:theo-harm-coho-corr-1}.
\end{thm}

In Section 5, the correspondence result is considered for non-Calabi-Yau case. Theorem \ref{sec5:theo-1} describes such a correspondence where only partial structure can be preserved.  

The first application of Theorem \ref{complete-corr} is to show the folklore result in physical literature that "the (topological field) theory on the relevant deformation part determines the theory on the marginal deformation part" is true for $tt^*$-geometry.

\begin{thm} Let $(\C^{n+2},f)$ be a LG model satisfying $\deg f=n+2$. Then the $tt^*$ geometrical structure on the relavant deformation space in $S_{mr}$  determines uniquely the pre-$tt^*$ structure $(H^{\CY}\to M, \kappa^{\CY},\eta^{\CY}, \D^{CY})$ on the Calabi-Yau side.
\end{thm}

\begin{proof} Since $M\subset S_{mr}$ is of higher codimension, the $tt^*$ bundle structure $\widehat{\E}_{\ominus}$ on $S_{mr}\setminus M$ can be extended uniquely to $M$ by taking the limits, in particular, the Higgs field $\widehat{C}$ can be extended holomorphicly to $M$.    
\end{proof}

For the convenience of the reader, we give a required description of the Gelfand Leray forms and the Gauss-Manin connections in Appendix \ref{Appendix-A}. 

\begin{rem} Theorem \ref{main-1} appeared as a part of the doctoral thesis of the second author \cite{L} based on the early work \cite{FLY}. In almost the same time, J. Yan and X. Tang \cite{TY} gave a different approach to this correspondence conjecture via the vaccum line bundles and Weil-petersson metrics. The approach and the result of this paper is different to \cite{TY}, which shows the correspondence from the $tt^*$ geometry of LG model to the (classical) $tt^*$ geometry. Our paper constructs a new $tt^*$ structure in CY side compared to the classical one. Since the $tt^*$ connections are only metric connections which are not unique a priori. It believes that the new $tt^*$ structure constructed in our paper is a deformation of the classical one. Another advantage of our approach is that one can easily build the correspondence for non-Calabi Yau hypersurface.

Another possible way to treat the LG and CY models is to use the polyvector fields and trace maps (see \cite{LLS}).     
\end{rem}

\subsection*{Acknowledgement} We would like to thank Junrong Yan and XinXing Tang for suggestions to our previous preprint kindly. The first author thanks Si Li and Emanuel Scheidegger for many useful discussions. The third author would like to thank Konstantin Aleshkin for helpful discussions.

\section{$tt^*$-structures of Landau-Ginzburg and Calabi-Yau models}

\subsection{Differential geometry of LG model}\

The LG model has been studied for a long time by physicists as an important model in topological field theories (\cite{Ce1,Ce2}). A systematic study of the differential geometrical structure of LG models appeared in \cite{Fan} by the first author. Let's recall the definitions and main results in \cite{Fan} and \cite[Appendix A]{FF}.

\begin{df}[{\cite[Definition 2.39]{Fan}}]\label{sec2-defn-LG} The LG model $(M,g,f)$ of dimension $N$ consists of a complex $N$-dimensional manifold $(M,g)$ and a superpotential function $f$ satisfying the conditions:
\begin{enumerate}
\item $(M,g)$ is a non-compact complete Kaehler manifold with metric $g$ having bounded geometry and
\item $f$ is a nontrivial holomorphic function on $M$.
\end{enumerate}
The LG model is said to be strongly tame, if for any
constant $C>0$, there is
\begin{equation}
|\nabla f|^2-C|\nabla^2 f|\to \infty, \;\text{as}\;d(x, x_0)\to
\infty.
\end{equation}
Here $d(x,x_0)$ is the distance between the point $x$ and the base
point $x_0$.
\end{df}

\begin{rem} Definition \ref{sec2-defn-LG} can be understood as the Kahler version of the LG model. One can generalize it to the \textquotedblleft complex \textquotedblright LG model without metric involved. In this case, one can study the complex deformation theory (\cite{KKP}). On the other hand, one can study the LG model with the action of a symmetry group. In \cite[Definition 2.39]{Fan}, LG models are called the section-bundle systems and more "tame" conditions has been discussed. The \textquotedblleft strongly tame condition \textquotedblright here was called as \textquotedblleft elliptic condition \textquotedblright in \cite{KL} and such form appeared much earlier in the study of 1-dimensional Schrodinger equations.
\end{rem}

Compared to the Cauchy-Riemman operators $\bpat, \pat$ on a compact Kahler manifold, one can study the corresponding twisted operators on a LG model $(M,g,f)$ of dimension $N$:
$$
\bar{\partial}_f=\bar{\partial}+\partial f\wedge,\quad\partial_f=\partial+\bar{\partial} \bar{f}\wedge.
$$
The Hodge $\star$ operator is defined as a $\C$-linear operator $\star: \Lambda^{p,q}\to \Lambda^{n-q,n-p}$ such that
$$
g(\varphi, \psi) d \text{vol}_M=\varphi\wedge \star \bar{\psi}.
$$
Then the $L^2$-conjugate operator $\bpat_f^\dag$ of $\bpat_f$ can be expressed as $\bpat_f^\dag=-\star \pat_{-f}\star$. Similarly, the $L^2$-conjugate operators of $\pat_f$ is $\pat_f^\dag=-\star \bpat_{-f}\star$. We have the twisted Laplacian
$$
\Delta_f=\bar{\partial}^{\dag}_f\bar{\partial}_f+\bar{\partial}_f\bar{\partial}^{\dag}_f.
$$
The commutativity of the Hodge $\star$ operators with the twisted operators has been carefully analyzed in \cite{Fan}. An important observation is the following identity found in \cite{Fan}:
$$
\star \Delta_f=\Delta_{-f}\star.
$$
This property has been used in \cite{FF} to prove the vanishing of the first Zeta function related to $\Delta_f$. We list some properties of the twisted operators as follows.

\begin{prop}[{\cite[Chapter 2]{Fan}}]\label{commutator} Let $(M,g,f)$ be a LG model, then we have
$$\partial_f^2=\bar{\partial}_f^2=0, \quad \bar{\partial_f}\partial_f+\partial_f\bar{\partial}_f=0,$$
$$(\bar{\partial}_f^{\dag})^2=(\partial^{\dag}_f)^2=0, \quad\bar{\partial}^{\dag}_f\partial^{\dag}_f+\partial^{\dag}_f\bar{\partial}^{\dag}_f=0,$$
$$[\partial_f,\bar{\partial}^{\dag}_f]=[\bar{\partial}_f,\partial^{\dag}_f]=0,$$
$$[\partial_f,\partial_f^{\dag}]=[\bar{\partial}_f,\bar{\partial}^{\dag}_f]=\Delta_f,$$
and the Kahler-Hodge identities:
$$[\partial_f,\Lambda]=-i\bar{\partial}^{\dag}_f, \quad [\bar{\partial}_f,\Lambda]=i\partial^{\dag}_f,$$
$$[\partial^{\dag}_f,L]=-i\bar{\partial}_f, \quad [\bar{\partial}^{\dag}_f,L]=i\partial_f,$$
where L is the Lefschetz operator and $\Lambda$=$*^{-1}\circ$L$\circ*$.  The twisted Laplacian $\Delta_f$ commutes with all the above operators.
\end{prop}

The spectrum of the twisted Laplacian associated to a strongly tame LG model has the following nice property (Note that the proof of Theorems \ref{sec2-theo-2.40}, \ref{sec2-theo-2.52} and \ref{sec2-theo-2.66} can also be found in \cite[Appendix A]{FF}.):

\begin{thm}[{\cite[Theorem 2.40]{Fan}}]\label{sec2-theo-2.40} Let $(M,g,f)$ be a strongly tame LG model of dimension $N$. Then $\Delta_f$ has purely discrete spectrum and all the eigenforms form a complete basis of the Hilbert space $L^2(\Lambda^{*}(M))$.
\end{thm}

Let $\mathcal{H}_f\subseteq$Dom($\Delta_f$) be the subspace of $\Delta_f$-harmonic forms, $E_\mu$ be the eigenspace with respect to the eigenvalue $\mu$ and $\Pi_\mu$ be the projection from $L^2(\Lambda^{*}(M))$ to $E_{\mu}$.  We have the spectrum decomposition formulas:
$$L^2(\Lambda^{*}(M))=\mathcal{H}_f\oplus(\oplus_{i=1}^\infty E_{\mu_i}),\quad \Delta_f=\sum_i\mu_i \Pi_{\mu_i}.$$
The Green operator $G_f$ of $\Delta_f$ satisfies:
$$
G_f\Delta_f+\Pi=\Delta_fG_f+\Pi=I,
$$
where we set $\Pi=\Pi_0$.

This implies the following Hodge-de-Rham decomposition:

\begin{thm}[{\cite[Theorem 2.52]{Fan}}]\label{sec2-theo-2.52} There are orthogonal decompositions for any $k=0,1,\dots,2N$:
$$
L^2\Lambda^k=\mathcal{H}^k_f\oplus im(\bar{\partial}_f)\oplus im(\bar{\partial}^{\dag}_f).
$$
In particular, we have the isomorphism
$$
H^*_{((2),\bar{\partial}_f)}\cong \mathcal{H}^*_f,
$$
where $H^*_{((2),\bar{\partial}_f)}$ is the cohomology of the following $L^2$-complex:
$$
\cdots\rightarrow L^2\Lambda^{k-1}\xrightarrow{\bar{\partial}_f}L^2\Lambda^k\xrightarrow{\bar{\partial}_f}L^2\Lambda^{k+1}\rightarrow\cdots.
$$
\end{thm}

The $L^2$-cohomology is given by the following result.

\begin{thm}[{\cite[Theorem 2.66]{Fan}}]\label{sec2-theo-2.66} Let $(M,g,f)$ be a strongly tame LG model of dimension $N$ and assume that $M$ is a Stein manifold. Then

\begin{equation}
\dim \mathcal{H}^k_f=
\begin{cases}
0,\;k\neq N\\
\mu,\;k=N.
\end{cases}
\end{equation}
and there is an explicit isomorphism:
\begin{equation}
\I:\mathcal{H}^N_f \to \Omega^N (M)/df\wedge\Omega^{N-1}(M),
\end{equation}
where $\Omega^*(M)$ are spaces of holomorphic forms on $M$.
\end{thm}

\subsubsection*{\underline{The construction of $\I$ and the induced real structure $\kappa_f$ on $R_f$}}\

Let $\rho$ be a smooth function with compact support in $M$ which equals to $1$ in a neighborhood of $Crit(f)$, the set of critical points of $f$. Define the following operator
$$
V_f= \sum^n_{i=1} \frac{\bar{f_i}}{|\nabla f|^2} (dz_i \wedge)^* : \quad \Omega^{*,*}(M \setminus Crit(f)) \rightarrow \Omega^{*-1,*}(M \setminus Crit(f)).
$$
A direct calculation gives the following result (or see \cite[Lemma A3 and A4]{FF}).
\begin{lm}
\begin{equation}
 [df \wedge, V_f]=1
 \end{equation}
 and
 \begin{equation}
 [\bar{\partial},[\bar{\partial},V_f]]=[df \wedge,[\bar{\partial},V_f]]=[V_f,[\bar{\partial},V_f]]=0.
 \end{equation}
Let
\begin{align*}
T_{\rho}= \rho + (\bar{\partial} \rho)V_f \frac{1}{1+[\bar{\partial},V_f]},\;
R_{\rho}= (1-\rho)V_f \frac{1}{1+[\bar{\partial},V_f]},
\end{align*}
then we have
\begin{equation}\label{equa-expl-cons-i}
[\bar{\partial}_f,R_{\rho}]=1-T_{\rho} \quad on \quad \Omega^*(M).
\end{equation}
\end{lm}
For any $[a]\in R_f$, we have a holomoprhic $N$-form $[A]=[adz_1\wedge\cdots\wedge dz_N]$. By (\ref{equa-expl-cons-i}), we have 
\begin{equation}
T_\rho A=A+\bpat_f (-R_\rho A). 
\end{equation}
Since $T_\rho A$ is a smooth compactly supported $N$-form, it has a unique harmonic $N$-form $\alpha_A$ representing the $L^2$-cohomological class $[T_\rho A]$. Now the following equation has a unique solution $\beta_{A,\rho}$ in the domain of $\bpat_f$:
\begin{equation}
\begin{cases}
\bpat_f \beta_{A,\rho}=\alpha_A-T_\rho A\\
\bpat^\dag_f\beta_{A,\rho}=0.
\end{cases}
\end{equation}
$\alpha_A$ has the holomorphic representative:
$$
\alpha_A=A+\bpat_f \eta_A,
$$
where the $N-1$ form $\eta_A=-R_\rho A+\beta_{A,\rho}$ has polynomial growth. Note that $A$ in the above representative is unique up to a term in $df\wedge \Omega^{N-1}(M)$. We define the map $\I_-([A])=\alpha_A$.  

On the other hand, to define $\I(\alpha)$ for any harmonic $N$-form $\alpha$ we need use the following identity on pairings. 

\begin{thm}[{\cite[Theorem 3.4]{FS}}]\label{thm-pair-FS}
Suppose $f\in\mathbb{C}[z_1, \dots, z_{N}]$ is a nondegenerate quasi-homogeneous polynomial, $\alpha, \beta$ are $\Delta_f$-harmonic $N$-forms on $\mathbb{C}^N$. Then there are polynomials $A, B\in\mathbb{C}[z_1, \dots, z_{N}]$ and $(N-1)$-forms $\mu, \nu$ such that $\alpha=Adz_1\wedge\dots\wedge dz_{N}+\bar{\partial}_f\mu$ and $\beta=Bdz_1\wedge\dots\wedge dz_{N}+\bar{\partial}_f\nu$. We have the following identity:
$$
\eta(\alpha,\beta):=\int_{\mathbb{C}^N}\alpha\wedge*\beta=k_N{\rm Res}_f(AB),
$$
where  $k_N=\frac{(-1)^{N(N-1)/2}i^N}{2^N}$.
\end{thm}
Note that in the above theorem, $\Res_f:=(2\pi i)^N \res_{f,0}$, where $\res_{f}$ is the residue appearing in the complex analysis (ref. \cite[Appendix A]{FLY}).

Choose a $\C$-basis $\{a_1,\cdots,a_\mu\}$ of $R_f$, we have a family of holomorphic $N$-forms $\{A_1,\cdots, A_\mu\}$, where $A_i=a_i dz_1\wedge\cdots\wedge dz_N$ and harmonic forms $\{\alpha_1,\cdots,\alpha_\mu\}$. Let $\eta_{ij}=\eta(\alpha_i, \alpha_j)$ and $(\eta^{ij})$ be the inverse matrix. For any harmonic $N$-form $\alpha$, we define 
\begin{equation}
\I(\alpha)=\sum_{i,j}\eta(\alpha,\alpha_i)\eta^{ij}A_j. 
\end{equation}
Using Theorem \ref{thm-pair-FS}, it is easy to prove that $\I\circ \I_-=\I_-\circ \I=I$ and $\I_-=\I^{-1}$. We call such a basis $\{\alpha_i\}$ corresponding to $\{a_i\}$ a holomorphic basis.  

Let $\hat{\kappa}^{\LG}$ be the usual complex conjugate, then it induces a real structure $\kappa_f$ acting on $R_f$ via the isomorphism $\I$. We define $\kappa_f$ as follows. Let $\{\alpha_i\}$ be a holomorphic basis. Then for any $k$, we have 
$$
\I(\hat{\kappa}^{\LG}(\alpha_k))=\eta(\overline{\alpha_k},\alpha_i)\eta^{ij}A_j.
$$

Denote by $\eta_{\kb i}:=\eta(\overline{\alpha_k},\alpha_i)$ and 
\begin{equation}\label{sec2:real-stru-1}
\K^j_{\kb}=\sum_i \eta_{\kb i}\eta^{ij}.
\end{equation} 
Then the action of $\kappa_f$ is defined as 
\begin{equation}
\kappa_f(A_k)=\sum_j \K^j_{\kb} A_j, \;\text{and}\;\kappa_f(\sum_k \lambda_k A_k)=\sum_k \overline{\lambda_k}\K^j_{\kb} A_j,
\end{equation}
for any complex numbers $\lambda_1,\dots,\lambda_k$. 

Let $\tilde{\alpha}_\mu$ be another holomorphic basis corresponding to the holomorphic $N$-form $\tilde{A}_\mu$. Denote by $\phi_{\mu i}=\eta(\alpha_\mu,\alpha_i)$, and $\kappa_f(\tilde{A}_\mu)=\sum_\nu \tilde{\K}^\nu_{\bar{\mu}}\tilde{A}_\nu$. Then we have the transformation formula
\begin{equation}
\tilde{\K}^\nu_\mub=\sum_{\ib, j}\phi_{\mub \ib}\K^{\ib }_j\phi^{\nu j},
\end{equation}
where 
$$
\phi_{\mub \ib}:=\overline{\phi_{\mu,i}},\; \K^{\ib}_j=\sum_{\mb, l}\eta^{\ib \mb}\K^l_{\mb}\eta_{lj}=\overline{\K^i_{\bar{j}}},\;\eta^{\ib \mb}:=\overline{\eta^{im}}.
$$
It is easy to prove the following conclusions:
\begin{lm}
\begin{enumerate} We have
\item \begin{equation}
\K^{\jb}_k \K^l_{\jb}=\delta^l_k,\;\text{or}\;\overline{\K}\cdot \K=I.
\end{equation}
\item $\alpha_A$ is a real harmonic $N$-form if and only if $\K(A)=A$. 
\item Let $A_i\in R_f$ be a family of holomorphic $N$-forms corresponding to a basis consisting of real harmonic $N$-forms $\alpha_i$, then the matrix $\K$ representing $\kappa_f$ in terms of the base $\{A_i\}$ is the identity matrix.  
\end{enumerate}
\end{lm}

\subsubsection*{\underline{The LG model $(\C^{N}, f)$}}\

We will apply the above conclusions to a family of LG models $(\C^{N}, f_u)$, where $\C^{N}$ is the Euclidean space and $f_u$ is a family of quasi-homogeneous polynomials.

Let $f:\C^{N}\to \C$ be a quasi-homogeneous polynomial with weights $(q_1,\cdots,q_N)$ if for all $\lambda\in \C^*$, there is
\begin{equation}
f(\lambda^{q_1}z_1,\cdots,\lambda^{q_N}z_N)=\lambda f(z_1,\cdots,z_N).
\end{equation}

\begin{df}\label{df-nondeg-quasi}
Let $f\in\mathbb{C}[z_1, \dots, z_N]$ be a quasi-homogeneous polynomial, it is called non-degenerate if\\
(1) $f$ contains no monomial of the form $z_iz_j$ for $i\neq j$.\\
(2) $f$ has only an isolated singularity at the origin.
\end{df}

It was shown in \cite[Proposition 2.1.6]{FJR} that for non-degenerate quasi-homogeneous polynomial, each variable $z_i$ has weight $q_i\le 1/2$.

The universal unfolding of $f$ can be described by the Milnor ring:
$$
R_f:=\mathbb{C}[z_1,...,z_{N}]/I_f,
$$
where $I_f=\left \langle \frac{\partial f}{\partial z_1}, \dots, \frac{\partial f}{\partial z_N} \right \rangle$ is the ideal generated by the derivatives of $f$. When $f$ is a non-degenerate quasi-homogeneous polynomial, the Milnor ring $R_f$ is finite-dimensional and its dimension $\mu_f=\dim R_f$ is called the Milnor number. Let $\{\phi_1,...,\phi_{\mu}\}$ be a basis of $R_f$ consisting of monomials,   and consider the following deformation of $f$:
$$
F(z, u)=f(z)+\sum_{j=1}^{\mu}u_j\phi_j(z).
$$
We denote by $u=(u_1, \dots, u_{\mu})$ the deformation parameter. The above deformation gives the universal deformation of $f$ and we have a family of LG models $(\C^N, F(z,u))$ defined on the deformation space $S\subset \C^\mu$, where $0\in S$.

\begin{df}
The deformation parameters $u_j$ are divided into three types by the weights of $\{\phi_j(z)\}$. $u_j$ is called:\\
(1) relevant, if the weight of $\phi_j$ is positive;\\
(2) marginal, if the weight of $\phi_j$ is zero;\\
(3) irrelevant, if the weight of $\phi_j$ is negative.\\
If a deformation direction is a linear combination of more than one deformation parameters,  we take the highest weight to be the weight of this direction.
\end{df}

The strong deformation of a general strongly tame LG models has been considered in \cite[Section 3.1.2]{Fan}. Here we only give the definition of the strong deformation of the LG model $(\C^N, f)$, where $f$ is a holomorphic function. Let $F(z,u), u\in S$, be a deformation of $f$ over the deformation space $S$ and denote $f_u(z)=F(z,u)$.

\begin{df}\label{df-strong-deformation} A family of LG models $(\C^N, f_u)$ over $S$ is called a strong deformation of $(\C^N, f)$, if $f_0=f(z)$, and the following conditions hold:
\begin{enumerate}
\item for any $u\in S$, $(\C^N, f_u)$ is a strongly tame LG model,
\item $\text{sup}_{u\in S}\mu(f_u)<\infty$,
\item for any $u\in S$, $\Delta_{f_u}$ have common domains in the space of $L^2$ forms.
\end{enumerate}
\end{df}

\begin{thm}[{\cite[Theorem 2.43, Section 3.1.3]{Fan}}]\label{sec2:theo-stro-defo-1}  Any LG model $(\C^N, f)$ with $f$ being a non-degenerate quasi-homogeneous polynomial is a strongly tame LG model. The marginal and relevant deformations of $f$ are strong deformations. Hence for all such LG models, Theorem \ref{sec2-theo-2.40}, \ref{sec2-theo-2.52} and \ref{sec2-theo-2.66} hold. If let $S_{mr}$ represent the parameter space consisting of the marginal and relevant deformations, then for any $u\in S_{mr}$ there is
$$
\dim \mathcal{H}_{f_u}^k=
\begin{cases}
0, &  k \neq N\\
\mu_f, & k=N.
\end{cases}
$$
and there exists an explicit isomorphism:
$$
\I:\mathcal{H}_{f_u}^N\longrightarrow \Omega^N(\mathbb{C}^N)/d f_u\wedge\Omega^{N-1}(\mathbb{C}^N)\cong R_{f_u}.
$$
Here $R_{f_u}$ is the Milnor ring of $f_u$.
\end{thm}

\subsection{$tt^*$-structures on LG models}\

Let $f(z)\in\mathbb{C}[z_1, \dots, z_{N}]$ be a non-degernerate quasi-homogeneous polynomial and consider the strong deformation $F(z,u)$ given in Theorem \ref{sec2:theo-stro-defo-1}:
\begin{equation}\label{sec2:equa-stro-1}
F(z, u)=f(z)+\sum_{i=1}^{s}u_i\phi_i(z), u \in S_{mr}.
\end{equation}
We assume that $S_{mr}\ni 0$ and has dimension $s$.

Theorem 2.4 of \cite{Fan} gives a $tt^*$ structure over $S_{mr}$ with respect to $f_u$. After applying the construction to the deformation $f_u/2$ and normalizing the pairing in this $tt^*$ structure, we get the so called big $tt^*$ structure in this paper.

\begin{thm}[Big $tt^*$ structure]\label{big-structure}\

Let $(\C^N, f_u)$ over $S_{mr}$ be the strong deformation of $(\C^N, f)$ given by (\ref{sec2:equa-stro-1}). Then there exists a $tt^*$ structure over $M$ denoted by
$$
\widehat{\E}^{\LG}=(\hat{H}^{\LG}\to M, \hat{\kappa}^{\LG},\hat{\eta}^{\LG}, \hat{D}^{\LG},\hat{\bar{D}}^{\LG}, \hat{C}^{\LG}, \hat{\bar{C}}^{\LG}).
$$
These data are given as follows:
\begin{enumerate}
\item $S_{mr}$ is the parameter space of the strong deformation $F$ consisting of the relevant and marginal directions.
\item $\hat{H}^{\LG}$ is a smooth complex vector bundle over $M$, and at $u\in S_{mr}$ the fiber $\hat{H}^{\LG}_u=\mathcal{H}^{N}_{f_u/2}$ consists of the $\Delta_{\frac{f_u}{2}}$-harmonic $N$-forms. $\hat{H}^{\LG}$ is called the Hodge bundle over $S_{mr}$.
\item $\hat{\kappa}^{\LG}$ is the real structure given by the usual complex conjuagate.
\item $\hat{\eta}^{\LG}$ is the pairing on a fiber $\hat{H}_u^{\LG}$ defined by
$$
\hat{\eta}^{\LG}(\alpha,\beta)(u)=\frac{1}{i^{(N-2)^2}(2\pi i)^N \mu}\int_{\mathbb{C}^N}\alpha(u)\wedge* \beta(u),
$$
where $\mu$ is the Milnor number of $f$ and $\alpha,\beta$ are two sections of $\hat{H}^{\LG}$. The $tt^*$ metric is given by
$$
\hat{g}^{\LG}(\alpha,\beta)=\hat{\eta}^{\LG}(\alpha,\hat{\kappa}^{\LG}\beta).
$$
\item The $tt^*$ connections are defined as
$$
\hat{D}_i^{\LG}=\Pi\circ\partial_i, \quad\hat{\bar{D}}_{\bar{i}}^{\LG}=\Pi\circ\bar{\partial}_{\bar{i}}, \quad i=1,\cdot\cdot\cdot,s,
$$
where $\Pi: L^2(\Lambda^N(\mathbb{C}^N))\rightarrow \mathcal{H}^N_{f_u/2}$ is the projection.
\item The Higgs fields are defined as
$$
\hat{C}_i^{\LG}=\Pi\circ\partial_iF,\quad\hat{\bar{C}}_{\bar{i}}^{\LG}=\Pi\circ\overline{\partial_iF}, \quad i=1,\cdot\cdot\cdot,s.
$$
\end{enumerate}
\end{thm}

\begin{rem}
 Since any $\Delta_{f_u}$-harmonic form is exponentially decaying at infinity (ref. \cite[Theorem 3.4.3]{Fan}), the pairing $\hat{\eta}^{\LG}(\alpha,\beta)(u)=\int_{\mathbb{C}^n}\alpha_u\wedge* \beta_u$ between two $\Delta_{f_u}$-harmonic forms $\alpha_u,\beta_u$ is well-defined. The smoothness of the bundle $\hat{H}^{\LG}$ is the conclusion of the stability theorem, \cite[Theorem 3.53]{Fan}.
\end{rem}

We want to give an explicit formula to the pairing. Firstly by Theorem \ref{thm-pair-FS}, we have the following conclusion.

\begin{crl}\label{new-LG-residue}
For any $\Delta_{\frac{f}{2}}$-harmonic form $\alpha$ and $\beta$ with the representation $\alpha=Adz_1\wedge\cdots\wedge dz_N+\bar{\partial}_{\frac{f}{2}}\mu$ and $\beta=B dz_1\wedge\cdots\wedge dz_N+\bar{\partial}_{\frac{f}{2}}\nu$, we have
\begin{equation}
\int_{\mathbb{C}^{N}}\alpha\wedge *\beta=i^{(N-2)^2}{\Res_{f}}(A B).
\end{equation}
Hence, we have
\begin{equation}
\hat{\eta}^{\LG}(\alpha,\beta)(u)=\frac{1}{\mu}\res_{f_u,0}( \I(\alpha)\I(\beta)).
\end{equation}
\end{crl}

There is a special pair $[1]$ and $[1^\vee]=[\det(\frac{\pat^2 f}{\pat z_i \pat_{z_j}})]$ in the Milnor ring $R_f$ such that $$
\Res_f([1][1^\vee])=(2\pi i)^N \mu,
$$
and
\begin{equation}
\hat{\eta}(\I^{-1}([1]), \I^{-1}([1^\vee]))=i^{(N-2)^2}(2\pi i)^N \mu.
\end{equation}

We remark that the pairing $\hat{\eta}^{\LG}$ defined here is a normalization to the original one in the $tt^*$ structure of \cite[Theorem 2.4]{Fan}.

\subsection{Induced $tt^*$ structure on cohomology bundles}\

\

Let $f_u(z)=F(z,u)$ be the deformation defined in (\ref{sec2:equa-stro-1}). Given $\alpha > 0$,  we can define two sets ${f_u}^{\geq \alpha}$ and ${f_u}^{\leq -\alpha}$:
\begin{displaymath}
{f_u}^{\geq \alpha}=\{z\in \mathbb{C}^{N}\lvert Re{f_u}(z)\geq \alpha\},\quad {f_u}^{\leq -\alpha}=\{z\in \mathbb{C}^{N}\lvert Re{f_u}(z)\leq -\alpha\}
\end{displaymath}
Since $F$ is a strongly deformation, for any $\alpha,\beta > 0$, the two sets ${f_u}^{\ge \alpha}$ and ${f_u}^{\ge \beta}$ are homotopic, and similarly for ${f_u}^{\le -\alpha}$ and ${f_u}^{\le -\beta}$. We define these two homotopic equivalent classes by ${f_u}^{\geq +\infty}$ and ${f_u}^{\leq -\infty}$.

By the exact sequence of homology groups
\begin{displaymath}
\cdots\rightarrow H_k({f_u}^{\geq +\infty},\mathbb{Z})\rightarrow H_k(\mathbb{C}^{N},\mathbb{Z})\rightarrow H_k(\mathbb{C}^{N},{f_u}^{\geq +\infty},\mathbb{Z})\rightarrow  H_{k-1}({f_u}^{\geq +\infty},\mathbb{Z})\rightarrow\cdots,
\end{displaymath}
we get the isomorphism
$$
H_k(\mathbb{C}^{N},{f_u}^{\geq +\infty},\mathbb{Z})\cong H_{k-1}({f_u}^{\geq +\infty},\mathbb{Z}).
$$
Similarly we have
$$
H_k(\mathbb{C}^{N},{f_u}^{\le -\infty},\mathbb{Z})\cong H_{k-1}({f_u}^{\le -\infty},\mathbb{Z}).
$$
It is known that the $N$-th homology group is the only non-vanishing homology group.

\begin{df} Define $\hch_{\ominus,u}:=H^N(\C^N,{f_u}^{-\le \infty}, \R)$ and $\hch_{\oplus,u}:=H^N(\C^N,{f_u}^{\ge \infty}, \R)$. Let
$\hch_{\ominus},\hch_{\oplus}$ be the bundles over $S_{mr}$ with fiber at $u\in S_{mr}$ to be $\hch_{\ominus,u}$ and $\hch_{\oplus,u}$ respectively. 
\end{df}

Let $\{\alpha_i,i=1,\cdots,\mu\}$ be a local frame of $\hat{H}^{\LG}$ consisting of the $\bpat_{F/2}$-harmonic $N$-forms. Then $\alpha_i$ are primitive forms which satisfy
$$
\bpat_{f_u/2} \alpha_i(u)=0,\;\pat_{f_u/2} \alpha_i(u)=0, \forall u\in S_{mr}. 
$$
\begin{lm}(ref. \cite[Lemma 4.47, Theorem 4.60]{Fan}) Let $S_i^-(u)=e^{(f_u+\overline{f_u})/2}\alpha_i$ and
$S_i^+(u)=e^{-(f_u+\overline{f_u})/2}*\alpha_i$. Then $S_i^-(u)$ and $S_i^+(u)$ are
$d$-closed $N$-forms on $\C^N$, furthermore, $\{S_i^-\}$ (or $\{S_i^+\}$) forms a flat frame of $\hch_{\ominus}$ (or $\hch_{\oplus}$) with respect to the topological Gauss-Manin connection $\D^{\top}$.  
\end{lm}

\begin{proof} We can prove by a direct calculation that $S_i^-(u)$ and $S_i^+(u)$ are closed. The topological Gauss Manin connection $\D^\top$ for the bundle $\hch_{\ominus}$ is given by the action on the basis $S_j^-(u),j=1,\dots,\mu$, and we have for any $\tau=1,\dots,\dim S_{mr}$,
\begin{align*}
\D^\top_\tau S_j^-(u)&=\frac{\pat}{\pat u_\tau} S_j^-(u)=e^{\frac{f_u+\overline{f_u}}{2}}(\pat_\tau \alpha_j+\pat_\tau (f_u/2) \alpha_j)\\ 
&=e^{\frac{f_u+\overline{f_u}}{2}}\left((\hat{D}_\tau^{\LG}+\hat{C}_\tau^{\LG})\alpha_j+(\bpat_{f_u/2}+\pat_{f_u/2})\bpat_{f_u/2}^*G(\pat_\tau (f_u/2))\alpha_j \right)\\
&=d\left(e^{\frac{f_u+\overline{f_u}}{2}}\bpat_{f_u/2}^*G(\pat_\tau(f_u/2)))\alpha_j\right),
\end{align*}
i. e. $\D^\top_\tau [S_j^-(u)]=0$.
\end{proof}
Due to the stability theorem, \cite[Theorem 3.53]{Fan}, the cohomological bundles $\hch_{\ominus}$ and $\hch_{\oplus}$ admit smooth structures.  

There is a natural isomorphism 
$$
\hat{\star}: \hch_{\ominus,u}\to \hch_{\oplus,u}: e^{(f_u+\overline{f_u})/2}\alpha_j\mapsto e^{-(f_u+\overline{f_u})/2}\star \alpha_j
$$
such that $\hat{\star}^2=\star^2=(-1)^N$. There is a pairing $\eta^\top$ on $\hch_{\ominus,u}$ defined by
$$
\eta^\top(S_k^-(u),  S_j^-(u)):=\frac{1}{i^{(N-2)^2}(2\pi i)^N \mu}\int_{\C^N} S_k^-(u)\wedge \hat{\star}S_j^-(u)=\hat{\eta}^{\LG}(\alpha_k, \alpha_j).
$$

The above discussion implies the following conclusions (ref. \cite[Theorem 4.61]{Fan}).

\begin{thm}\label{sec2:theo-harm-coho-corr-1}
The map 
$$
\Phi: \alpha(u)\mapsto [e^{(f_u+\overline{f_u})/2}\alpha(u)]
$$ 
gives a $tt^*$ bundle isomorphism from $\hat{\E}=(\hat{H}^{\LG}, \hat{\eta}^{\LG}, \D, \hat{\kappa}^{\LG})$ to $\hat{\E}_{\ominus}:=(\hch_{\ominus},\eta^\top, \D^{\top},\hat{\kappa}^{\LG})$. Here $\hat{\kappa}^{\LG}$ is the usual complex conjugate. 
\end{thm}

\begin{proof} We only need to check that $\Phi\circ \D=\D^{\top}\circ \Phi$. This is easy to see since we can split $\D^{\top}=D^{\top}+C^{\top}+\bar{D}^{\top}+\bar{C}^{\top}$ by defining
$$ 
D_\tau^{\top} S_j^-(u):=e^{(f_u+\overline{f_u})/2} \hat{D}_\tau^{\LG}\alpha_j,\; C^{\top}_\tau S_j^-(u):=e^{(f_u+\overline{f_u})/2} \hat{C}_\tau^{\LG}\alpha_j,
$$
where $D^{\top},\bar{D}^{\top}$ are $tt^*$ connections of $\eta^{\top}(\cdot, \hat{\kappa}\cdot)$. 
\end{proof}

\subsubsection*{\underline{Gauge transformation and holomorphic basis}}\

Let $\{\alpha_j(u),j=1,\dots,\mu\}$ be a holomorphic basis of $\hat{H}_u^{\LG}$ such that $I(\alpha_j (u))=A_j=a_j dz_1\wedge\cdots\wedge dz_N$, where $a_j\in R_f$. We define 
$$
C_{i j}^k(u)=\eta_u(a_i A_j, A_l)\eta_u^{lk}
$$ 
which depends holomorphically on $u$. Then the action of $C^{\top}_\tau$ is complex linear which is given by
$$
C^{\top}_\tau\cdot S_j^-(u)=C_{\tau j}^k(u)  S_k^-(u),\forall \tau=1,\dots,s, \forall j=1,\dots,\mu.
$$
Here $s=\dim S_{mr}$.

Denote by $C=\sum_{\tau=1}^{s} C_\tau d\tau$, where $C_\tau=(C_{\tau j}^k)$ is a $\mu\times \mu$ matrix. $C$ is a holomorphic section of $\Omega^{1,0}(S_{mr}, \End(\hch_{\ominus}))$. Similar to \cite[Proposition 4.63]{Fan}, we have the following conclusion.
\begin{lm}\label{sec2:lemm-gaug-tran-1} There is a smooth section $\Theta$ of the bundle $\End( \hch_{\ominus})$ satisfying the following equation:
\begin{equation}\label{sec2-gaug-tran-1}
\begin{cases} 
\pat\Theta=C\\
\bpat\Theta=0,
\end{cases}
\end{equation}
with $\Theta(0)=0$. Moreover, the following identity holds:
\begin{equation}\label{sec2-gaug-tran-2}
[C_\tau, \Theta]\equiv [\bar{C}_\tau, \Theta]\equiv 0, \forall \tau=1,\dots,m.
\end{equation}
\end{lm}

\begin{proof} Note that $D^{\top}=\pat-C^{\top}$ and $C^{\top}$ satisfy the $tt^*$ equations $D^{\top}\circ D^{\top}=0$ and $C^{\top}\wedge C^{\top}=0$.  We have the following equalities: 
$$
\pat C\equiv 0, \;[C,C]\equiv 0. 
$$
This showes that the equation (\ref{sec2-gaug-tran-1}) is integrable and by applying the $\bpat$-Poincare lemma there exists a unique solution. Since 
$$
\pat(C\Theta-\Theta C)=[\pat C,\Theta]-[C,C]\equiv 0,
$$
and $\Theta(0)=0$, we can prove $[C_\tau, \Theta]\equiv 0$. Using the facts that $[\bar{D}^{\top}, C]=0$ and $\bpat C=0$, we can prove $[\bar{C}_\tau, \Theta]\equiv 0$. Then (\ref{sec2-gaug-tran-2}) is proved.
\end{proof}

By Lemma \ref{sec2:lemm-gaug-tran-1}, we have the following result. 
\begin{prop}\label{sec2: prop-gaug-tran-1}
$S$ is a holomorphic local section of $\hch_{\ominus}$ satisfying $D^{\top}S=0$ if and only if $e^{-\Theta(u)}S$ is a holomorphic flat section of $\D^{\top}$. 
\end{prop}
\begin{proof} The conclusion is due to the following gauge transformation:
$$
e^{-\Theta(u)}\circ D^{top}\circ e^{\Theta(u)}=e^{-\Theta(u)}\circ (\pat-C)\circ e^{\Theta(u)}=\pat.
$$
Note that we also have 
\begin{equation}
e^{-\Theta(u)}\circ \bar{D}^{top}\circ e^{\Theta(u)}=\bpat-\bar{C}=\bar{D}^{top}.
\end{equation}
\end{proof}

\begin{crl}\label{sec2:crol-gaug-tran-1} Take any basis $\{A_j,j=1,\dots,\mu\}$ in the Milnor ring $R_f$. Let $\alpha_j(u)=\alpha_{A_j}(u)$ be the corresponding $\Delta_{f_u/2}$-harmonic forms in $\hat{H}^{\LG}_u$. The following conclusions hold:
\begin{enumerate}
\item $S^-(u)=(S^-_1(u),\dots,S^-_\mu(u))^T$ is a local flat frame of $\hch_{\ominus}$, where 
$$
S^-_j(u)=e^{\frac{f_u+\overline{f_u}}{2}}\alpha_j(u),
$$ 
and $e^{\Theta(u)}S^-(u)$ is a holomorphic frame of $\hch_{\ominus}$ and is horitontal w. r. t. $D^{\top}$. 
\item Define
$$
S_{\ominus,j}(u)=e^{f_u/2}A_j, j=1,\dots,\mu,  
$$
and $S_{\ominus}(u)=(S_{\ominus,1},\dots, S_{\ominus,\mu})^T$. Then $S_{\ominus}(u)$ is a holomorphic frame of $\hch_{\ominus}$ which is horizontal w. r. t. $D^{\top}$ and 
$e^{-\Theta(u)}S_{\ominus}(u)$ is a local flat frame of $\hch_{\ominus}$.
\item There is a constant matrix $\Xi$ such that $[S^-(u)]=\Xi\cdot e^{-\Theta(u)}\cdot [S_{\ominus}(u)]$. Here 
\begin{equation}
\Xi=(\int_{\Gamma^-_{{1}^\vee}}S_{\ominus}(0),\cdots, \int_{\Gamma^-_{{\mu}^\vee}}S_{\ominus}(0)),
\end{equation}
where $\PD(\Gamma^-_{{k}^\vee})\in \hch_{\oplus,0}$ is a dual basis of $S^-_j(0)$ such that
$$
\int_{\C^N}\PD(\Gamma^-_{{k}^\vee})\wedge S^-_j(0)=\delta_{kj}.
$$
\item The real structure $\kappa^{\LG}$ has the following matrix representation in terms of the holomorphic basis $\{S_{\ominus,j}, j=1,\dots,\mu\}$:
\begin{equation}\label{sec2:crol-gaug-tran-eq-4}
\K^{\top}=\overline{\Xi}\cdot e^{-\overline{\Theta(u)}}\cdot \K\cdot e^{\Theta(u)}\cdot \Xi^{-1}.
\end{equation}
\end{enumerate}
\end{crl}

\subsubsection*{\underline{Calculation of the period integrals}}\

In this part $f:\mathbb{C}^{N}\longrightarrow \mathbb{C}$ is assumed to be a nondegenerate homogeneous polynomial of degree $d$.

\begin{lm}[{\cite[Lemma 2.5]{FLY}}]\label{Lef-thim}
There exists a basis of $H_{N}(\mathbb{C}^{N},f^{\leq -\infty},\mathbb{Z})$ consisting of $\mu$ Lefschetz thimbles $\{\Gamma^{-}_a,a=1,...,\mu\}$ such that the images of these Lefschetz thimbles under $f$ are the negative real line.  More specifically, taking $\mu$ vanishing spheres $S_a(z)$ on $V_{-1}=f^{-1}(-1)$, then $\Gamma^{-}_a=\{S_a(t^{\frac{1}{d}}z)| t\ge 0\}$ are such Lefschetz thimbles.
\end{lm}

\begin{rem} Note that $\Gamma_a^-$ is invariant under the scaling deformation $z\mapsto T^{\frac{1}{d}}z$ for any $T>0$, hence $\Gamma^{-}_a$ are also the Lefschetz thimbles of $H_{N}(\mathbb{C}^{N},{(\frac{f}{T})}^{\leq -\infty},\mathbb{Z})$.
\end{rem}

\begin{lm}\label{harmonic-holomorphic-integral}
Let $\{\Gamma^{-}_a,a=1,...,\mu\}$ be the Lefschetz thimbles given in Lemma \ref{Lef-thim}, and $A=\beta dz_1\wedge \cdots dz_{N}$, where $\beta\in R_f$ is a homogeneous element. Then
$$
\int_{\Gamma^{-}_a}e^{f}A=\Gamma(\frac{N+\deg[\beta]}{N})
\int_{S_a}\frac{A}{df},
$$
where $S_a=\Gamma^-_a \cap f^{-1}(-1)$ is the vanishing sphere. 
\end{lm}

\begin{proof} This is a direct application of the Gelfand-Leray form associated to a holomorphic $N$-form and the scaling transformation of the Lefschetz thimbles. For the definition and the properties of the Gelfand-Leray form, the reader can refer to \cite{AGV} or \cite[Appendix C]{FLY}.
\end{proof}

\begin{rem}\label{coef-discuss}
There is a mistake appeared in \cite[Theorem 2.6]{FLY} which claims that the integration of $S^-_j$ is equal to the integration of $S_{\ominus,j}$ along the Lefschetz thimble. This mistake comes from \cite[Lemma 4.88]{Fan} because $e^{f+\bar{f}}\bpat_f R$ is not a $d$-closed form and one can't use Poincare duality. We thank J. Yan and X. Tang to point out this.The correct statement is in Corollary \ref{sec2:crol-gaug-tran-1}.

Compared to the results in \cite{FLY}, we use $f/2$ instead of $f$ to construct the Hodge bundle over $S_{mr}$. This scaling can simplify the coefficient appearing in the integral of the Gelfand-Leray form.
\end{rem}

\subsection{Small $tt^*$ structure on the LG side}\

In this section,  we assume that $N=n+2$ and $f=f(z_1,\dots,z_{n+2}):\mathbb{C}^{n+2}\longrightarrow\mathbb{C}$ is a non-degenerate homogeneous polynomial with degree $n+2$. The Milnor ring $R_f$ is graded by the polynomial degree and is of $\mu$-dimension as a $\C$-vector space. We denote by $R_f^{k}$ the degree $k$ part and the subring $R_f^{(n+2)*}=\oplus_{i=0}^{\infty} R_f^{(n+2)i}$.

In \cite{FLY}, a $tt^*$ substructure of the big $tt^*$ structure of the LG model was constructed and is called the small $tt^*$ structure. We recollect the procedures to build the small $tt^*$ structure.

The following theorem in \cite{FLY} is proved by combining Lemma \ref{harmonic-holomorphic-integral} with a calculation of Gauss-Manin connection.
\begin{thm}[{\bf Theorem} 2.14 \cite{FLY}]
Assume $f:\mathbb{C}^{n+2}\longrightarrow\mathbb{C}$ is a non-degenerate homogeneous polynomial of degree $n+2$. Then the cohomology classes
$$
\{\frac{a dz_1\wedge \cdots dz_{n+2}}{df}|\; a \in R^{(n+2)*}_f\}
$$
are exactly the invariant parts of $H^{n+1}(V_{-1},\mathbb{C})$ under the action of the monodromy group of Gauss-Manin connection. Moreover, these classes span a real subspace of $H^{n+1}(V_{-1},\mathbb{C})$.
\end{thm}

The following result is a slight modification of \cite[Theorem 2.9]{FLY}, with $f$ replaced by $f/2$.

\begin{defthm}[Small $tt^*$ structure]\label{thm-small-lg} Let $M\subset S_{mr}$ be the marginal deformation part.
Define a bundle $H^{\LG}$ over $M$, where the fiber at $u\in M$ is given by
$$
H^{\LG}_u=\{\alpha\in \mathcal{H}^{n+2}_{f_u/2}|\; i(\alpha)\in R_{f_u}^{(n+2)*}\}.
$$
Then there is natural embedding of bundles:
$$
\begin{diagram}
H^{\LG} & \rTo & \hat{H}^{\LG}\\
\dTo & & \dTo\\
M &\rTo &S_{mr}.
\end{diagram}
$$
The actions $\hat{\kappa}^{\LG}, \hat{D}^{\LG},\hat{\bar{D}}^{\LG}, \hat{C}^{\LG}, \hat{\bar{C}}^{\LG}$ from the big $tt^*$ structure preserve the subbundle $H^{\LG}$. Hence we get a $tt^*$ substructure
$$
\E^{\LG}=({H}^{\LG}\to M, {\kappa}^{\LG},{\eta}^{\LG}, {D}^{\LG}, {\bar{D}}^{\LG}, {C}^{\LG}, {\bar{C}}^{\LG}),
$$
which we call it the small $tt^*$ structure on the LG side.
\end{defthm}

Correspondingly, via the $tt^*$ bundle isomorphism $\Phi: \widehat{\E}^{\LG}\to \widehat{\E}_{\ominus}$, the small $tt^*$ subbundle $\E^{\LG}$ is mapped to the $tt^*$ subbundle 
\begin{equation}
\E_{\ominus}=(\cho\to M, {\kappa}^{\LG}, \eta^{\top}, D^{\top}, \bar{D}^{\top}, C^{\top},\bar{C}^{\top}).
\end{equation} 

Denote
\begin{equation}
\dim M=\dim R_f^{n+1}=m,\quad \dim R_f^{(n+2)*}=\mu_s.
\end{equation}

\subsection{A $tt^*$ structure on CY side}\

In this section, we will formulate the (original) $tt^*$ structure appeared in the literatures (ref. \cite{St}, or \cite[Theorems 2.11-2.14]{FLY}). This $tt^*$ structure was proved to have partial correspondence to the small $tt^*$ structure of the LG model (\cite{FLY}). However, we will built later a new $tt^*$ structure based on this original one which exactly coincides to the small $tt^*$ structure in LG side under the big residue map which we will define these objects later on.

Let $f=f(z_1,\dots, z_{n+2}):\C^{n+2}\rightarrow\C$ be a non-degenerate homogeneous polynomial with $\deg f=n+2$. Then the hypersurface $X_f\subseteq\mathbb{CP}^{n+1}$ defined by $f=0$ is a Calabi-Yau $n$-fold. The deformation space of $X_f$ in $\mathbb{CP}^{n+1}$ can be identified with the marginal deformation space $M$ of the LG model $(\C^{n+2},f)$. Note that if $n\neq 2$, $M$ can be identified with the deformation space of complex structures on $X_f$, and if $n=2$, $\dim M$ is less than $1$ of the deformation space as a complex manifold (ref. \cite[Section 3.1.4]{Fan}).

There is a filtration on the primitive Hodge bundle $H^{\CY}$ over $M$ whose fiber at $u\in M$ is $H^{\CY}_u=H^n_{prim}(X_{f_u}, \C)$:

$$
F^n (H^{\CY})^n\subseteq\dots \subseteq F^0 (H^{\CY})^n=H^{\CY},\quad F^k (H^{\CY})^n=\oplus_{p=0}^{n-k}(H^{\CY})^{n-p, p}.
$$
and the complex conjugation induces a real structure on $H^{\CY}$.

Since $H^{\CY}$ is the tensor product of a local system with $\C$ over $M$, it induces the Gauss-Manin connection $\D^{\CY}$. Denote by $\nabla^{\CY}$ and $\overline{\nabla}^{\CY}$ the $(1,0)$-part and $(0,1)$-part of $\D^{\CY}$.

For any $ p\in \N, \alpha\in (H^{\CY})^{p,n-p}$, there are the decompositions:
\begin{align*}
\widetilde{D}(\alpha)&=\Pi_{(H^{\CY})^{p,n-p}}[\nabla^{\CY}(\alpha)],\;\overline{\widetilde{D}}(\alpha)=\Pi_{(H^{\CY})^{p,n-p}}[\overline{\nabla}^{\CY}(\alpha)],\\
\widetilde{C}(\alpha)&=\Pi_{(H^{\CY})^{p-1,n-p+1}}[\nabla^{\CY}(\alpha)],\;\overline{\widetilde{C}}(\alpha)=\Pi_{(H^{\CY})^{p+1,n-p-1}}[\overline{\nabla}^{\CY}(\alpha)].
\end{align*}
where $\Pi_V$ are the orthogonal projection operators to the subspace $V$.

By the Griffiths transversality theorem, we have the decompositions:
$$
\nabla^{\CY}=\widetilde{D}+\widetilde{C},\quad\overline{\nabla}^{\CY}=\overline{\widetilde{D}}+\overline{\widetilde{C}}.
$$
The flatness of $\D^{\CY}$ implies Cecotti-Vafa's $tt^*$-equations. By Hodge-Riemann bilinear relation, the hermitian inner product
$$
g_p(u, v)=i^{2p-n}\int_Xu\wedge\bar{v}
$$
gives a hermitian metric on $H^{n-p, p}_{prim}$. 

Hence $g=\sum_p g_p$ is a hermitian metric on $H_{prim}^n$ and $\widetilde{\eta}:=g(\cdot,\kappa^{\CY}\cdot)$ is a nondegenerate pairing.

We have the following well-known result.

\begin{thm}[ref. {\cite[Theorem 2.12]{FLY}}]\label{tt*-CY}
Let $f:\C^{n+2}\rightarrow \mathbb{C}$ be a non-degenerate polynomial of degree $n+2$. Then the data
$$
\widetilde{\E}^{\CY}=(H^{\CY}\to M, \kappa^{\CY},\widetilde{\eta}, \widetilde{D}, \overline{\widetilde{D}},\widetilde{C}, \overline{\widetilde{C}})
$$
forms a  $tt^*$ structure.
\end{thm}

\subsection{Small residue map and Frobenius algebras correspondence}\

Carlson and Griffith \cite{CG} studied the residue map, and used it to construct an explicit isomorphism between a subring of the Milnor ring and the $tt^*$ bundle on CY side. There is a well-known topological interpretation of the residue map which we recall as follows.

For any ($k-1$)-cycle $\gamma$ on $X_f$, the image $T(\gamma)$ under the Leray coboundary map $T:H_{k-1}(X_f)\longrightarrow H_k(\mathbb{CP}^{n+1}-X_f)$ is the boundary of an $\epsilon$-tubular neighborhood of $\gamma$, for a small $\epsilon>0$. Let $\alpha\in H^k(\mathbb{CP}^{n+1}-X_f)$. The (topological) residue $\Res$ is defined as the formal adjoint of $T$:
$$
\int_{\gamma}\Res\alpha=\int_{T(\gamma)}\alpha,
$$
where $\Res(\alpha)$ is a $(k-1)$-dimensional cohomology class defined on $X_f$. The topological residue $\Res$ is $2\pi i$ times the analytical residue.

Define a holomorphic form
$$
\Omega=\sum_{i=1}^{n+1}(-1)^iz_idz_1\wedge\dots\wedge\widehat{dz_i}\wedge\dots\wedge dz_{n+2}
$$
on $\mathbb{C}^{n+2}$. Note that the $\Omega$ defined as above is different to the $\Omega$ considered in \cite{FLY} by an opposite sign.

For any holomorphic polynomial $A$, the rational form $\Omega_A=\frac{A\Omega}{f^{a+1}}$ is a meromorphic $(n+1)$-form with $X_f$ as its polar locus. When the degree of $A$ is chosen to make the quotient homogeneous of degree zero, i.e. $\deg A=(n+2)a$, then the rational form $\Omega_A$ is said to have adjoint level $a$.

\begin{thm}[{\cite[Chapter 3]{CG}}]\label{prop-residual-map}

Let $\Omega^{n+1}(pX_f)$ be the sheaf consisting of meromorphic $(n+1)$-forms with at most poles of order $p$ on $X_f$ and let $\Gamma \Omega^{n+1}(pX_f)$ denote the set of its global sections. Then we have
\begin{enumerate}
\item $\Res \Gamma\Omega^{n+1}((n+1)X_f)=H^n_{prim}(X_f,\mathbb{C})$.
\item $\Res\Gamma\Omega^{n+1}((a+1)X_f)=F^{n-a}H^n_{prim}(X_f,\mathbb{C})$.
\item let $\Omega_A$ be of adjoint level $a$, then
$$
r(A):=\Res\Omega_{A}
$$ has Hodge level $n-a+1$ if and only if $A$ lies in the Jacobian ideal of $f$.
\end{enumerate}
\end{thm}

\begin{df}[Small residue map]\label{eq:def-of-residual-4}
We call the map
$$
r: R_f^{(n+2)*}\to H_{prim}^n (X_f,\C)
$$
the small residue map.
\end{df}

\begin{crl}[{\cite[Theorem 3.2]{FLY}}] The map $A\mapsto (\Res\Omega_A)^{n-a,a}$ induces an isomorphism
$$
R_f^{(n+2)a}\longrightarrow H^{n-a,a}_{prim}(X_f),
$$
where $(\Res\Omega_A)^{n-a,a}$ denotes the $(n-a, a)$-part of $\Res\Omega_A\in H^n_{prim}(X_f,\mathbb{C})$.
\end{crl}

In \cite{FLY} we proved the following isomorphism between Frobenius algebras:

\begin{thm}[{\cite[Theorem 3.5-3.7]{FLY}}]\label{fro-iso}

Let $f\in\mathbb{C}[z_0, \dots, z_{n+1}]$ be a non-degenerate holomorphic homogeneous polynomial with degree $n+2$.

The $tt^*$ bundle $H^{LG}$ on LG side has a Frobenius algebra structure which is isomorphic to the Frobenious algebra on the Milnor ring $R_f$. The multiplication in $R_f$ is given by the polynomial multiplication $[A]\cdot[B]=[AB]$, the pairing is given by the residue pairing
$$
([A], [B]):=\Res_f(A, B),\quad {\rm deg}A+{\rm deg}B=(n+2)n.
$$
The $tt^*$ bundle $H^{\CY}$ on CY side has also a Frobenius algebra structure on $H^n_{prim}(X_f)=\oplus_{a=0}^nH^{n-a, a}(X_f)$. The multiplication is given by the Yukawa product
$$
[\alpha]*[\beta]:=\iota_\Omega({\iota_\Omega}^{-1}[\alpha]\wedge{\iota_\Omega}^{-1}[\beta]),
$$
where $\iota_\Omega$ is the contraction by the holomorphic volume form  $\Omega=(-1)^n(\Res 1)^{n, 0}$:
$$
\iota_\Omega: H^a(\wedge^a TX_f)\longrightarrow H^{n-a,a}(X_f),\quad [s]\mapsto [s\vdash\Omega],
$$
and the pairing is equivalent to the residual pairing under the bundle isomorphism $r$, i.e.  we have
$$
(\alpha, \beta):=\Res_f(A, B)=k_{ab}^{-1}\int_{X_f}r(A)\wedge r(B),
$$
which is nonzero only if $a+b=n$, and in this case we have
\begin{align*}
\alpha&=(\Res\Omega_A)^{n-a, a},\deg A=(n+2)a;\\
\beta&=(\Res\Omega_B)^{n-b, b},\deg B=(n+2)b.
\end{align*}
and
$$
k_{ab}=\frac{(-1)^{\frac{a(a+1)+b(b+1)}{2}+b^2+n}}{a!b!}(n+2).
$$
If we scale the above formula by multipling the constant $p_{n+2}=2^{n+2}i^{-n^2}$,  then the map below gives an isomorphism between Frobenius structures:
$$
r^{\prime}: \oplus_{a=0}^nR_f^{(n+2)a}\longrightarrow H^n_{prim}(X_f)=\oplus_{a=0}^n H_{prim}^{n-a, a}(X_f)
$$
\begin{equation}\label{sec2:defn-smal-resi-1}
r^{\prime}(A)=c_a^{-1}(\Res\Omega_A)^{n-a, a}, \;\forall A \;{with}\; {\rm deg}A=(n+2)a,\quad c_a=\frac {(-1)^{n+a(a+1)/2}} {a!}.
\end{equation}
For the multiplication, up to a sign (which depends on $a,b,n$),  we have $r^{\prime}(AB)=r^{\prime}(A)r^{\prime}(B)$.  The pairings are equal up to a power of $i$ (which depends on $a,b,n$).
\end{thm}

\section{The big residue map}

In this section we will introduce another kind of residue map defined by Steenbrink in \cite{Ste}. This residue map is defined on the whole Milnor ring, which has larger domain than the residue map defined in \cite{CG} (which is only defined on the subring $R_f^{(n+2)*}$). We call the residue map from \cite{Ste} as the big residue map, and the previous one as the small residue map.

We will prove several formulas between big residue map, small residue map, Lefschetz thimbles and oscillating integrals. These formulas are important tools for the construction of the complete $tt^*$ correspondence.

\subsection{Compactification of Milnor fibers}\label{subsec:def-of-cpt}\

Let $f(z_1,\dots,z_{n+2}):\mathbb{C}^{n+2}\longrightarrow \mathbb{C}$ be a non-degenerate homogeneous holomorphic polynomial on $\mathbb{C}^{n+2}$,  and ${\rm deg}f=d$.  By adding a new variable $z_{n+3}$ and a complex number $t\in\mathbb{C}$,  we can define a family of homogeneous polynomials:
\begin{equation}
F_t(z_1,\dots,z_{n+2},z_{n+3})=f(z_1,\dots,z_{n+2})-tz_{n+3}^d.
\end{equation}
By the non-degeneracy of $f$,  it is easy to see that $F_t$ has no nontrivial critical point as $t\neq 0$. As $t\neq 0$, we denote by $\bar{V}_t \subset \mathbb{CP}^{n+2}$ the smooth hypersurface defined by $F_t=0$. Let $\mathbb{CP}^{n+1}\subset \mathbb{CP}^{n+2}$ be the hyperplane defined by $\{z_{n+3}=0\}$.

We denote the intersection of the $(n+2)$-dimensional complex plane $\mathbb{CP}^{n+2}-\{z_{n+3}=0\}\cong \mathbb{C}^{n+2}$ with $\bar{V}_t$ to be $V_t$. Under the local coordinates $\{x_i=\frac{z_i}{z_{n+3}},i=1\dots,n+2\}$ of $\mathbb{C}^{n+2}\subseteq\mathbb{CP}^{n+2}$, $V_t$ is defined by the following equation:
\begin{equation}
f(x_1,\dots,x_{n+2})-t=0.
\end{equation}
Hence $V_t$ is just the Milnor fiber $f_t=f^{-1}(t)$.

Note that
\begin{equation}\label{relation1}
X_f=\bar{V}_t\cap\{z_{n+3}=0\}=\bar{V}_t-V_t, \forall t\in\mathbb{C}^*.
\end{equation}
We also denote by $V_{\infty}:=X_f\hookrightarrow \C P^{n+1}$.

By the non-degeneracy of $f$,  we have the following conclusion.

\begin{lm}
Let $f:\mathbb{C}^{n+2}\rightarrow \C$ be a non-degenerate homogeneous polynomial with $\deg f=d$.  Then for any $t\in\mathbb{C}^*$, $\bar{V}_t \subset \C P^{n+2}$ is a smooth hypersurface and intersects the hyperplane $\mathbb{CP}^{n+1}=\{z_{n+3}=0\}$ transversally.
\end{lm}

\subsection{Big residue map}\

Given a smooth complex manifold $X$ and a smooth hypersurface $Y$ on $X$,  we have the Leray coboundary map $\delta$:
\begin{equation}\label{def:leray-coboundary-map}
\delta:H_{m-1}(Y)\longrightarrow H_m(X-Y)
\end{equation}
such that for any $(m-1)$-cycle $\alpha$ on $Y$, $\delta([\alpha])$ can be obtained by taking the boundary of a tubular neighborhood of $\alpha$ in the normal bundle of $Y$.

The topological residue map $Res$ is defined as the dual map of $\delta$: take any homology class $\alpha\in H_{m-1}(Y)$ and any cohomology class $\omega\in H^m(X-Y)$, we have
\begin{equation}\label{def:topo-residual-map}
\int_{\delta(\alpha)}\omega=\int_{\alpha}Res(\omega)
\end{equation}

Now we can give the definition of the big residue map:
\begin{df}[Big residue map]\label{big-residual-map}
Let $f(z_1,\dots,z_{n+2}):\mathbb{C}^{n+2}\rightarrow \C$ be a non-degenerate homogeneous polynomial with $\deg f=d$. Take a monomial basis $\{A_i\}_{i=1}^{\mu}$ of the Milnor ring $R_f$, and define the map $h: \{A_i\}_{i=1}^{\mu}\rightarrow \mathbb{Q}$ by
\begin{equation}\label{def-of-h}
h(A_i)=\frac{{\rm deg}A_i+n+2}{d},
\end{equation}
and the fractional form
\begin{equation}\label{def-of-omega}
\omega_{A_i,t}=(f-t)^{[-h(A_i)]}A_idz_1\wedge\dots\wedge dz_{n+2},\qquad t\in\mathbb{C}^*.
\end{equation}
The big residue map $R_t: R_f \rightarrow H^{n+1}(V_t)$ is defined by
\begin{align}
R_t(A_i)=\Res(\omega_{A_i,t}),\label{def-of-big-residual-map}.
\end{align}
where $\Res: H^{n+2}(\C^{n+2}-V_t)\rightarrow H^{n+1}(V_t)$ is the residue map.
\end{df}

\begin{thm}[{\cite[Theorem A.1 of Appendix A]{Ste} or \cite{CIR}}]\label{cir-conclusion} Under the topological residue map $\Res$, the set $\{\omega_{A_i,t}\}_{i=1}^{\mu}$ is mapped to a basis of $H^{n+1}(V_t)$. Moreover, under the Deligne filtration (see \cite{CIR})
$$
0=\mathcal{W}_{n}\subseteq \mathcal{W}_{n+1}\subseteq \mathcal{W}_{n+2}=H^{n+1}(V_t)
$$
the subset $\{\omega_{A_i,t}|\; h(A_i)\not\in \Z \}$ is mapped to a basis of $\mathcal{W}_{n+1}$, and the subset $\{\omega_{A_i,t}|\; h(A_i)\in \Z\}$ is mapped to a basis of $\mathcal{W}_{n+2}/\mathcal{W}_{n+1}$ (after projection).
\end{thm}

\begin{crl}
The big residue map $R_t: R_f \to H^{n+1}(V_t)$ is an isomorphism.
\end{crl}

\subsection{The relation between the big and small residue maps}\

The domain of the big residue map $R_t$ is the whole Milnor ring $R_f$, while the domain of the small residue map $r$ is a subring of $R_f$. The map $h$ from (\ref{def-of-h}) can be viewed as a grading function, and the subring is characterized as the subset $\{a|\; a\in R_f, h(a)\in\Z \}$. These two residue maps also have different ranges. The range of the big residue map is the cohomology group $H^{n+1}(V_t)$,  while the range of the small residue map is $H_{prim}^n(V_{\infty})$.

In this section we will show that the small residue map can be factorized through the big residue map, and the subset $\{a|\; a\in R_f, h(a)\not\in\Z \}$ vanishes naturally under the composition maps so that only the subring $\{a|\; a\in R_f, h(a)\in\Z \}$ survives.

Given $t\in\mathbb{C}^*$,  we can define four different residue maps. The first one is
\begin{equation}\label{eq:def-of-residual-1}
R_t: H^{n+2}(\mathbb{C}^{n+2}-V_t)\longrightarrow H^{n+1}(V_t).
\end{equation}
$R_t$ is just the topological version of the big residue map.

By (\ref{relation1}), we have
$$
V_{\infty} = \bar{V}_t\cap\{z_{n+3}=0\} = \bar{V}_t-V_t,\; \forall t\in\mathbb{C}^*,
$$
the second residue map is
\begin{equation}\label{eq:def-of-residual-2}
r_{2,t}:H^{n+1}(V_t)\longrightarrow H^{n}(V_{\infty}),
\end{equation}

For the open manifold $\mathbb{CP}^{n+2}-\bar{V}_t$, we have:
\begin{equation}\label{eq:relation-of-cpt-2}
(\mathbb{C}^{n+2}-V_t)=(\mathbb{CP}^{n+2}-\bar{V}_t)-(\mathbb{CP}^{n+1}-V_{\infty})
\end{equation}

So we can define the third residue map for the pair $(\mathbb{C}^{n+2}-V_t,\mathbb{CP}^{n+1}-V_{\infty})$:
\begin{equation}\label{eq:def-of-residual-3}
r_{3,t}:H^{n+2}(\mathbb{C}^{n+2}-V_t)\longrightarrow H^{n+1}(\mathbb{CP}^{n+1}-V_{\infty}).
\end{equation}

Finally for the hypersurface $V_{\infty}\subset \mathbb{CP}^{n+1}$, we have the residue map:
\begin{equation}
r_{4}:H^{n+1}(\mathbb{CP}^{n+1}-V_{\infty})\longrightarrow H^n(V_{\infty})
\end{equation}
used to define the small residue map $r$ in Definition \ref{eq:def-of-residual-4}.

Denote by $\{\delta_{1,t},\delta_{2,t}, \delta_{3,t},\delta_4\}$ the Leray coboundary maps which we have used to define $\{R_t,r_{2,t},r_{3,t}, r_4\}$.

\begin{rem}
When $t=0$,  $V_0$ and $\bar{V}_{0}$ are singular, but we still have
$$
V_{\infty}=\bar{V}_0\cap\mathbb{CP}^{n+1}=\bar{V}_0-V_0.
$$
So the residue map $r_{3,0}$ is still well-defined.
\end{rem}

\begin{lm}\label{residual-commu} Let $t\in \C^*$. For any $\alpha\in H^k(\mathbb{C}^{n+2}-V_t)$, we have
\begin{equation}\label{eq:residual-commu}
r_{2,t}R_t(\alpha)=-r_{4}r_{3,t}(\alpha)
\end{equation}
\end{lm}

\begin{proof} Since the residue map is the dual of the Leray coboundary map, it suffices to show that
\begin{equation}\label{eq-resi-comm-2}
\delta_{1,t}\delta_{2,t}=-\delta_{3,t}\delta_{4}(\beta),\quad \forall\beta\in H_{n}(V_{\infty}).
\end{equation}
This follows from the fact that for any $t\in \C^*$, the smooth hypersurface $\bar{V}_t$ is transversal to the hyperplane $\C P^{n+1}$ in $\C P^{n+2}$. Therefore the normal bundle of $V_\infty$ in $\C P^{n+2}$ splits into the direct sum of the normal bundle of $V_\infty$ in $\C P^{n+1}$ and the normal bundle of $V_\infty$ in $\bar{V}_t$. Note that $\delta_{2,t}(\beta)$ equals the cup product of the lift of $\beta$ and the angular form of the corresponding normal bundle, it is obvious that (\ref{eq-resi-comm-2}) holds.
\end{proof}

Consider the holomorphic vector field $X$ on $\mathbb{C}^{n+2}$:
\begin{equation}
X(z_1,\dots,z_{n+2})=-z_1\frac{\partial}{\partial z_1}-\cdots-z_{n+2}\frac{\partial}{\partial z_{n+2}}.
\end{equation}
and the following $(n+1)$ form $\Omega$ by contracting the holomorphic volume with $X$:
\begin{equation}\label{eq:def-of-Omega}
\Omega=l_X(dz_1\wedge\cdots\wedge dz_{n+2})=\sum_{i=1}^{n+2}(-1)^iz_idz_1\wedge\cdots\wedge\widehat{dz_{i}}\wedge\cdots\wedge dz_{n+2}.
\end{equation}

\begin{lm}
Let $f(x_1,\dots,x_{n+2}):\mathbb{C}^{n+2}\longrightarrow\mathbb{C}$ be a non-degeneate homogeneous polynomial of degree $n+2$, $A$ is a homogeneous polynomial with $h(A)=a+1$(i.e.  $h(A)(n+2)={\rm deg}A+n+2$). Then the following conclusions hold:
\begin{enumerate}
\item The meromorphic $(n+1)$-form
$$
\Omega^{\mathbb{C}^{n+2}}_{A}:=\frac{A\Omega}{f^{a+1}}
$$
has pole along $V_0=\{(z_1,\dots,z_{n+2})\in\mathbb{C}^{n+2}\lvert f(z_1,\dots,z_{n+2})=0\}$ and can be pushed down to a meromorphic $n+1$-form $\Omega^{\mathbb{CP}^{n+1}}_A$ in $\C P^{n+1}$ with poles along $V_\infty$.
\item Let
$$
\omega_{A,0}(x_1,\dots,x_{n+2})=\frac{A d x_1\wedge\dots\wedge d x_{n+2}}{f^{a+1}}
$$
be the form defined in $\C^{n+2}\subset \C P^{n+2}$ with the local coordinates $x_i=z_i/z_{n+3}$, then we have
\begin{equation}\label{eq:new-Omega-A-by-residual-3}
r_{3,0}(\omega_{A,0})=\Omega_A^{\mathbb{CP}^{n+1}},
\end{equation}
where
\begin{equation}\label{def-of-small-residual-map}
r_{3,0}:H^{n+2}(\mathbb{C}^{n+2}-V_0)\longrightarrow H^{n+1}(\mathbb{CP}^{n+1}-V_{\infty})
\end{equation}
is the residue map defined before.
\end{enumerate}
\end{lm}

\begin{proof}

\begin{enumerate}
\item Note that the projective space $\C P^{n+1}$ is the quotient of $\C^{n+2}-\{0\}$ under the obvious $\C^*$-action. If $A$ is a homogeneous polynomial with degree $\deg A=(n+2)a$, then the meromorphic form $\Omega^{\mathbb{C}^{n+2}}_{A}$ is a $\C^*$-invariant form and can be push down to $\C P^{n+1}$. Hence $\Omega^{\mathbb{C}^{n+2}}_{A}$ can be viewed as a closed meromorphic $n+1$-form with pole $V_\infty$ in $\C P^{n+1}$.
\item In homogeneous coordinates $(z_1,\dots, z_{n+3})$, we have
\begin{align*}
&\omega_{A,0}(x_1,\dots,x_{n+2})=\omega_{A,0}(\frac{z_1}{z_{n+3}},\dots,\frac{z_{n+2}}{z_{n+3}})\\
=&\frac{A(z_1,\dots, z_{n+2})}{f^{a+1}(z_1,\dots,z_{n+2})}(dz_1-z_1\frac{d z_{n+3}}{z_{n+3}})\wedge\cdots \wedge (dz_{n+2}-z_{n+2}\frac{d z_{n+3}}{z_{n+3}}).
\end{align*}
Hence we have
$$
r_{3,0}(\omega_{A,0})=\Omega_A^{\mathbb{CP}^{n+1}}.
$$
\end{enumerate}
\end{proof}

Hence we have the following relation between the residue maps.

\begin{thm}\label{residue-limit} Let $A$ be a homogeneous polynomial with $h(A)=\frac{{\rm deg}A+n+2}{n+2}$. Let $\omega_{A,t}=(f-t)^{[-h(A)]}Adx_1\wedge\dots\wedge dx_{n+2}$ be a meromorphic $(n+2)$-form on $\C^{n+2}\subset \mathbb{CP}^{n+2}$. The following conclusions hold for any $t\in\mathbb{C}^*$:
\begin{enumerate}
\item If $h(A)\in \Z$, we have
\begin{equation}\label{small-residual-by-limit}
lim_{t\rightarrow 0}r_{4}r_{3,t}(\omega_{A,t})=r_4 r_{3,0}(\omega_{A,0}).
\end{equation}

\item If $h(A)\not\in \Z$, then
\begin{equation}
r_{3,t}(\omega_{A,t})=0.
\end{equation}
\end{enumerate}
\end{thm}

\begin{proof}Let $\Res:\mathbb{CP}^{n+1}-V_{\infty}\longrightarrow V_{\infty}$ be the small residue map.
\begin{enumerate}

\item For any homology class $\alpha\in H_{n}(V_{\infty})$, we have
\begin{equation}\label{eq:residual-limit}
\begin{aligned}
&lim_{t\to 0}\int_{\delta_{4}(\alpha)}r_{3,t}(\frac{Adz_1\wedge...\wedge dz_n}{(f-t)^{[h(A)]}})\\
=&\int_{\delta_{4}(\alpha)}r_{3,0}(\omega_{A,0}),
\end{aligned}
\end{equation}
since the image $r_{3,t}(\frac{Adz_1\wedge...\wedge dz_n}{(f-t)^{[h(A)]}})$ converges uniformly to $r_{3,0}(\frac{Adz_1\wedge...\wedge dz_n}{f^{[h(A)]}})$ on any compact set in $\mathbb{CP}^{n+1}-V_{\infty}$ as $t\to 0$.
This proved (1).
\item In homogeneous coordinates $(z_1,\dots, z_{n+3})$, we have
\begin{align*}
&\omega_{A,t}(x_1,\dots,x_{n+2})=\omega_{A,t}(\frac{z_1}{z_{n+3}},\dots,\frac{z_{n+2}}{z_{n+3}})\\
=&\frac{A(z_1,\dots, z_{n+2})}{(f(z_1,\dots,z_{n+2})-tz_{n+3}^d)^{[h(A)]}}z_{n+3}^{d[h(A)]-\deg A-n-2}(dz_1-z_1\frac{d z_{n+3}}{z_{n+3}})\wedge\cdots \wedge (dz_{n+2}-z_{n+2}\frac{d z_{n+3}}{z_{n+3}}).
\end{align*}
Since $h(A)\not\in \Z$, $d[h(A)]-\deg A-n-2\neq 0$. When taking the residue along $z_{n+3}=0$, we have
$$
r_{3,t}(\omega_{A,t})=0.
$$
Hence we proved (2).
\end{enumerate}
\end{proof}

\section{A complete correspondence between $tt^*$ structures}\label{Sec4}

In this section, we will focus on the construction of a complete LG/CY correspondence between $tt^*$ structures. Firstly we will do some preparations by proving some results concerning the residue maps, Lefschetz thimbles, Gauss-Manin connections and the residual pairings. Based on these results, we can construct a new $tt^*$ structure on the CY side. Comparing with the old $tt^*$ structure on CY side, the new one has the same $tt^*$ bundle, Gauss-Manin connection and pairing. However, there is a big difference between the old and the new one is the way to decompose the Gauss-Manin connection into the sum of the Higgs fields and the $tt^*$ connections. This implies that the Hodge bundle of the new $tt^*$ structure is the complex deformation of the old one. Finally, we will prove that the small residue map induces the isomorphism between the small $tt^*$ structure on the LG side and the new $tt^*$ structure on the CY side.

\subsection{Some preparations}\

In this section, we assume that $f=f(z_1,\dots, z_{n+2}):\mathbb{C}^{n+2}\rightarrow\mathbb{C}$ is a non-degenerate homogeneous polynomial with degree $d$.

For any $t\in\mathbb{C}^*$, $V_t=f^{-1}(t)$ is the fiber of the Milnor fibration and its homology group $H_{n+1}(V_t,\mathbb{Z})$ is a free abelian group of rank $\mu$. When $t$ varies in $\C^*$, it gives a local system over $\C^*$. We have the $\C$-vector bundle $\mathcal{H}\longrightarrow \mathbb{C}^*$ with the fiber $H_{n+1}(V_t,\mathbb{C})=H_{n+1}(V_t,\mathbb{Z})\bigotimes_{\mathbb{Z}}\mathbb{C}$. The local system equips $\mathcal{H}$ with a flat structure and the induced Gauss-Manin connection $\D^{GM}$. The parallel transition induces the action of the monodromy group.

Note that under the compactification, $V_t\subset \bar{V}_t\subset \C P^{n+2}$ and $V_\infty=\bar{V}_t-V_t$. We can construct some invariant flat sections of $\mathcal{H}$ from the cycles of $V_\infty$.

\begin{lm}\label{flat-section}\label{inv-sec}
For any integral $n$-cycle $\gamma$ in $V_\infty$, $s_\gamma(t):=[\delta_{2,t}(\gamma)]$ defines a flat section of $\mathcal{H} \to \C^*$, where
$\delta_{2,t}:H_{n}(V_{\infty},\mathbb{Z})\longrightarrow H_{n+1}(V_t,\mathbb{Z})$ is the Leray coboundary map. Moreover, for any $n$-cycle $\gamma$ in $V_\infty$, the section $s_\gamma$ is invariant under the action of the monodromy group.
\end{lm}

\begin{proof} Since the Milnor fibration is locally trivial, for any $t\in \C^*$, the fiber of $\mathcal{H}$ near $t$ can be identified with $\mathcal{H}_t=H_{n+1}(V_t,\mathbb{C})$. Hence $s_\gamma(t)$ is a locally constant section of $H_{n+1}(V_t,\mathbb{Z})$ and is invariant under the parallel transition of the Gauss-Manin connection. In particular, it is invariant under the action of the monodromy group.
\end{proof}

Lemma \ref{flat-section} gives a method to construct flat sections of the Gauss-Manin connection. In Lemma \ref{Lef-thim}, we have  considered the parallel transition of vanishing spheres and the parallel transitions of the $\mu$ vanishing cycles of $V_{-1}$ form the $\mu$ Lefschetz thimbles.

Let $\{\Gamma_i^{-},i=1,\dots,\mu\}$ be the set of Lefschetz thimbles constructed in Lemma \ref{Lef-thim}, which forms a basis of $H_{n+2}(\mathbb{C}^{n+2},f^{-\infty},\mathbb{Z})$. The intersection set $\{\Gamma_i^{-}\bigcap V_{t},i=1,\dots,\mu\}$ consists of the vanishing cycles and gives a basis of $H_{n+1}(V_t,\mathbb{Z})$ for $t\in \R_-$.

By Lemma \ref{flat-section}, for any $[\gamma]\in H_{n}(V_{\infty},\mathbb{C})$, $s_\gamma$ is a flat section of $\mathcal{H}$. So it can be written as a complex linear combination of Lefschetz thimbles $\{\Gamma_i^{-},i=1,\dots,\mu\}$ with $\mu$ complex constants $a_1,\dots,a_{\mu}$:
\begin{equation}\label{sec4:flat-conn-equa-1}
\delta_{2,t}(\gamma)=\sum_{i=1}^{\mu}a_i[\Gamma_i^{-}\cap V_{t}],\quad \forall t\in\mathbb{R}_{-}.
\end{equation}

\begin{lm}\label{residue-integral} Let $[\gamma]\in H_{n}(V_{\infty},\mathbb{C})$ be any homology class and expressed by (\ref{sec4:flat-conn-equa-1}). For any homogeneous polynomial $A$ satisfying $h(A)=\frac{\deg A+n+2}{d}\in\mathbb{Z}$, we define
\begin{equation}
Q_i(A, t)=\int_{\Gamma^-_i\cap V_t}R_t(\omega_{A,t}),
\end{equation}
and a family of meromorphic forms
$$
\{\omega_{A,t}=(f-t)^{-h(A)}Adz_1\wedge\dots\wedge dz_{n+2}\}
$$
over $\C^{n+2}$ depending on $t\in \R_-$ (or $\R_+$). Then for any $t\in \R_-$ (or $\R_+$), we have
$$
Q_i(A, t)=\int_{\Gamma_i^{-}\cap V_{t}}\frac{Adz_1\wedge\dots\wedge dz_{n+2}}{df}=(-t)^{h(A)-1}Q_i(A, -1),
$$
and the integral identity 
\begin{equation}
\int_{\delta_{2,t}(\gamma)}R_t(\omega_{A,t})=(-1)^{h(A)-1}\sum_{i=1}^\mu a_i Q_i(A,-1),
\end{equation}
which is a constant independent of $t$.
\end{lm}

\begin{proof} It suffices to consider the case $t\in\mathbb{R}_{-}$. By Lemma \ref{Lef-thim}, $Q_i(A,t)$ can be integrated out directly via the Gelfand-Leray form:
\begin{equation}\label{sec4:eq-Q-1}
Q_i(A, t)=\int_{\Gamma_i^{-}\bigcap V_{t}}\frac{Adz_1\wedge\dots\wedge dz_{n+2}}{df}=(-t)^{h(A)-1}Q_i(A,-1), \forall t\in\mathbb{R}_{-}.
\end{equation}

On the other hand, we have
\begin{equation}\label{eq:def-of-integral-Q}
Q_i(A, t)=\int_{\delta_{1,t}(\Gamma_i^{-}\bigcap V_{t})}\frac{Adz_1\wedge\dots\wedge dz_{n+2}}{f-t},
\end{equation}
where $\delta_{1,t}$ is the Leray coboundary map $\delta_{1,t}:V_t\rightarrow\mathbb{C}^{n+2}-V_t$.

Let $U_t$ be a small tubular neighborhood of $V_t$ in $\C^{n+2}$. Then for $t'$ approaching $t$ sufficiently small, the homology class  $\delta_{1,t'}(\Gamma_i^{-}\bigcap V_{t'})$ represents the same homology class in $\{H_{n+2}(\mathbb{C}^{n+2}-U_t,\mathbb{C})\}$. So it is valid to take the derivatives through the integration. We have  
\begin{displaymath}
\begin{aligned}\
\frac{d^k}{dt^{k}}Q_i(A,t)&=\frac{d^k}{dt^{k}}\int_{\delta_{1,t}(\Gamma_i^{-}\bigcap V_{t})}\frac{Adz_1\wedge\dots\wedge dz_{n+2}}{f-t}\\
&=\int_{\delta_{1,t}(\Gamma_i^{-}\bigcap V_{t})}\frac{d^k}{dt^{k}}\frac{Adz_1\wedge\dots\wedge dz_{n+2}}{f-t}\\
&=k!\int_{\delta_{1,t}(\Gamma_i^{-}\bigcap V_{t})}\frac{Adz_1\wedge\dots\wedge dz_{n+2}}{(f-t)^{k+1}}.
\end{aligned}
\end{displaymath}
Hence by the definition of $R_t:H^{n+2}(\mathbb{C}^{n+2}-V_t,\mathbb{C})\longrightarrow H^{n+1}(V_t,\mathbb{C})$,  we have
\begin{align*}
\int_{\Gamma_i^{-}\bigcap V_{t}}R_t(\omega_{A,t})&=\int_{\delta_{1,t}(\Gamma_i^{-}\bigcap V_{t})}\omega_{A,t}=\frac{1}{(h(A)-1)!}\frac{d^{h(A)-1}}{dt^{h(A)-1}}Q_i(A,t)\\
=&(-1)^{(h(A)-1)} Q_i(A,-1).
\end{align*}
Finally, we have 
\begin{equation}
\int_{\delta_{2,t}(\gamma)}R_t(\omega_{A,t})=(-1)^{h(A)-1}\sum_{i=1}^\mu a_i Q_i(A,-1).
\end{equation}

\end{proof}

\subsection{Correspondence between bundles}\label{sec4-4.2}

From now on, we assume that $\deg f=d=n+2$.

We first construct the correspondence between the base spaces.

By Theorem 3.21 in \cite{Fan}, the complex dimension $m$ of the marginal deformation space and the dimension of the space of deformation space of complex structures on $X_f$ as a projective variety are the same and
$$
m={\rm dim}H^1(X_f,T_{X_f})=\dbinom{n+1+n+2}{n+2}-(n+2)^2.
$$.

The marginal deformation of the polynomial $f$ can be written as:
$$
F(z, u)=f(z)+\sum_{i=1}^mu_i\phi_i(z),\quad {\rm deg}(\phi_i)=n+2,\quad i=1, \dots, m.
$$

Assume that the deformation parameter $u=(u_1,\dots,u_m)$ varies in a small neighborhood $M$ of the origin $0\in\mathbb{C}^m$. For any $u \in M$, $X_{F(z,u)}$ is a smooth Calabi-Yau hypersurface in $\mathbb{CP}^{n+1}$. So we get a fibration  $Y\rightarrow M$ whose fiber at $u$ is $X_{F(z,u)}$.

We can identify a marginal polynomial $\phi_i$ with the tangent vector $\pat_{u_i}$ on $M$, and can be also identified with the infinitesimal deformation of complex structures induced by the deformation $\{X_{F(z,u_i\phi_i)}\}$.

The base space $M$ can also be identified with the $\C$-vector space $R_f^{n+2}$ in the Milnor ring $R_f$.

Now we consider the fiber correspondence. On the CY side, the fiber of the $tt^*$ bundle $H^{\CY}$ is given by the primitive cohomology $H_{prim}^n(X_{f_u},\mathbb{C})$. On the LG side, the fiber is given by the space of the harmonic $(n+2)$-form of the operator $\Delta_{\frac{f_u}{2}}$ corresponding to the subring $R_f^{(n+2)*}$. Comparing to the setting in \cite{FLY}, using the operator $\Delta_{\frac{f_u}{2}}$ instead of $\Delta_{{f_u}}$ can simplify our correspondence (see Theorem \ref{harmonic-holomorphic-integral} and Remark \ref{coef-discuss}).
%The difference of the two settings are listed below:

%\begin{itemize}
%\item A polynomial which is a holomorphic representation of a $\Delta_{f_u}$-harmonic form may not be a holomorphic representation of a
%$\Delta_{\frac{f_u}{2}}$-harmonic form.

%\item By Remark \ref{coef-discuss}, the constant $\frac{\Gamma(\frac{n+2+deg[i(\alpha_a)]}{m})}{2^{\frac{n+2+deg[i(\alpha_a)]}{m}}}$ in Theorem \ref{harmonic-holomorphic-integral} changes to a simpler expression $\Gamma(\frac{n+2+deg[i(\alpha_a)]}{m})$ if we replace $f$ by $\frac{f}{2}$.
%\end{itemize}

We construct the correspondence between fibers by applying the small residue map. We denote by $V_{\infty,u}\subset \C P^{n+1}$ the hypersurface defined by the zero locus of $f_u(z)$.

Recall that $\delta_{4,u}: H_{n}(V_{\infty,u},\mathbb{C})\rightarrow H_{n+1}(\mathbb{CP}^{n+1}-\mathbb{V}_{\infty,u},\mathbb{C})$ is the Leray coboundary map dual to the residue map $r_{4,u}:H^{n+1}(\mathbb{CP}^{n+1}-\mathbb{V}_{\infty,u},\mathbb{C})\rightarrow H^{n}(V_{\infty,u},\mathbb{C})$ defined as:
$$
\int_{\gamma}r_{4,u}(\alpha)=\int_{\delta_{4,u}(\gamma)}\alpha,
$$
for any $\gamma\in H_n(V_{\infty,u},\mathbb{C}), \alpha\in H^{n+1}(\mathbb{CP}^{n+1}-V_{\infty,u},\mathbb{C})$.

For a homogeneous polymonial $A\in R_{f_u}$ with $\deg A=(n+2)a$, we have a meromorphic $(n+1)$-form $\frac{A\Omega}{f_u^{a+1}}$ on $\mathbb{C}^{n+2}$, where
$$
\Omega=\sum_{i=1}^{n+2}(-1)^{i}z_idz_1\wedge\dots\wedge\hat{dz_i}\wedge\dots\wedge dz_{n+2}.
$$
This form can be pushed down to a closed $(n+1)$-form $\Omega_{A,u}^{\mathbb{CP}^{n+1}}$ on $\mathbb{CP}^{n+1}-V_{\infty,u}$. Hence by applying the small residue map, we obtain a cohomology class
$$
r(A)=r_{4,u}(\Omega_A^{\mathbb{CP}^{n+1}}) \in H^n_{prim}(X_{f_u},\mathbb{C}).
$$

The fiber of the $tt^*$ bundle in the small $tt^*$ structure on the LG side is isomorphic to the subring $R_{f_u}^{(n+2)*}$ in the Milnor ring $R_{f_u}$. In the new $tt^*$ structure on CY side, we still choose the fiber of the $tt^*$ bundle $H^{\CY}$ to be the primitive cohomology of $V_{\infty,u} \subset \mathbb{CP}^{n+1}$. By Theorem \ref{prop-residual-map}, the map $r: R_{f_u}^{(n+2)*}\rightarrow H^n_{prim}(X_{f_u},\mathbb{C})$ is an isomorphism, and the image $\{r(A_i)=r_{4,u}(\Omega_{A_i}^{\mathbb{CP}^{n+1}}),i=1,\dots, \mu_s\}$ forms a basis of $H^n_{prim}(X_{f_u},\mathbb{C})$.

We can construct a local trivialization of the $tt^*$ bundles near $u=0$ by taking $A_i$ to be generators of $R_f$ before taking the small residue maps.

To give a better explict bundle isomorphism, we modify the definition of the small residue map $r$ slightly by scaling.

\begin{df}\label{def-residue} The modified residue map $\Rb:R_{f_u}^{(n+2)*}\rightarrow H^n_{prim}(X_{f_u},\mathbb{C})$ is defined as: for any $A_i\in R_{f_u}^{(n+2)*}$, there is
\begin{equation}\label{eq:def-of-r-prime-prime}
\Rb(A_i)=(-1)^aa!\Res(\Omega_{A_i}^{\mathbb{CP}^{n+1}}),\quad {\rm deg}A_i=(n+2)a.
\end{equation}

Instead of constructing an isomorphism between $\E^{\LG}$ between $\E^{\CY}$ directly, we construct an isomorphism between $\Eo$ and $\E^{\CY}$. 

The correspondence map $\Rbb: \cho \to \mathcal{H}$ is given by the composition of the following bundle isomorphism at any $ u\in M$:
\begin{equation}
\ch_{\ominus,u}\xrightarrow{\Io} R_{f_u}^{(n+2)*}\xrightarrow{\Rb} H^{\CY}_u,
\end{equation}
\end{df}
where $\Io^{-1}(A):=[e^{f_u}A]$.
\begin{df} Let $\{A_j,j=1,\dots,\mu_s\}$ be a basis of $R_f^{(n+2)*}$, and define 
\begin{equation}
S_{\ominus,j}(u)=\Io^{-1}(A_j),\;S^{\CY}_j(u)=\Rb(A_j).
\end{equation}
We call $S_{\ominus}(u)=(S_{\ominus_1}(u),\dots,S_{\ominus_{\mu_s}}(u))^T$ and $S^{\CY}(u)=(S^{\LG}_1(u),\dots,S^{\CY}_{\mu_s}(u))^T$ the holomorphic frames of the bundles $\cho$ and $H^{\CY}$. 
\end{df}

\subsection{New $tt^*$ bundle structure on CY side}\label{sec4-4.4}\

\

We observe that: for any $[\gamma(u)]\in H_n(X_{f_u},\mathbb{Z})$, there is
\begin{equation}\label{eq:cal-11}
\begin{aligned}\
&\int_{\gamma(u)}\partial_{u_\tau}\Rb(A_j)=(-1)^aa![\partial_{u_\tau}\int_{\delta_3(\gamma)}\frac{A_j\Omega}{F^{a+1}(u,z)}]\\
=&(-1)^aa![\int_{\delta_3(\gamma)}\partial_{u_\tau}(\frac{A_j\Omega}{(f+\sum_{k=1}^{m}u_k\phi_k)^{a+1}})]=(-1)^{a+1}(a+1)![\int_{\delta_3(\gamma)}\frac{\phi_\tau A_j\Omega}{(f+\sum_{k=1}^{m}u_k\phi_k)^{a+2}}]\\
=&(-1)^{a+1}(a+1)![\int_{\gamma}\Res\Omega^{\mathbb{CP}^{n+1}}_{\phi_\tau A_j,u}]=\int_{\gamma}\Rb(\phi_\tau A_j)=\int_{\gamma}S^{CY}_{\phi_\tau A_j}
\end{aligned}
\end{equation}
Here $\delta_3=\delta_{3,0}: H_{n+1}(\mathbb{CP}^{n+1}-V_{\infty,u}) \to H_{n+2}(\mathbb{C}^{n+2}-V_{0,u})$ is the Leray coboundary map and we have used the identity (\ref{eq:new-Omega-A-by-residual-3}). Hence we have 
$$
\partial_{u_\tau}\Rb(A_j)=S^{CY}_{\phi_\tau A_j}=\Rb(\phi_\tau A_j)=\sum_{k=1}^{\mu_s} C_{\tau j}^k S^{\CY}_k.
$$
\begin{df}[New $tt^*$-structure in CY side]\label{def-connection-higgs} 
We define the Gauss-Manian covariant derivative $\D^{\CY}_\tau=\pat_\tau$ and the Higgs field $C^{\CY}_{\tau}$ by their actions on the holomorphic sections $S^{\CY}_j (u)$:
$$
C^{\CY}_\tau\cdot S^{\CY}_j (u)=\sum_{k=1}^{\mu_s} C_{\tau j}^k S^{\CY}_k=\Rb(\phi_\tau A_j).
$$ 
The $tt^*$ connection $D^{\CY}$ is defined by $D^{\CY}_\tau=\pat_\tau-C_\tau$. The real structure $\kappa^{\CY}$ is the anti-complex linear map  defined by the action on the holomorphic basis:
\begin{equation}
\kappa^{\CY}\cdot S^{\CY}_j(u)=\sum_{k=1}^{\mu_s} (\K^{\top})^k_{\jb}(u) S^{\CY}_k(u),
\end{equation}
where the matrix $\K^{\top}$ is defined in (\ref{sec2:crol-gaug-tran-eq-4}). Via the real structure $\kappa^{\CY}$, we can define the $(0,1)$-parts of $\D^{CY}$ as follows:
\begin{equation}
\bar{D}^{\CY}=\kappa^{\CY}\circ D^{\CY}\circ \kappa^{\CY},\; \bar{C}^{\CY}=\kappa^{\CY}\circ C \circ \kappa^{\CY}.
\end{equation}
Now we have the new $tt^*$-structure $\E^{\CY}=(H^{\CY}\to M,\kappa^{\CY}, D^{\CY},\bar{D}^{\CY}, C^{\CY}, \bar{C}^{CY})$.
\end{df}

Analogous to Corollary \ref{sec2:crol-gaug-tran-1}, we have 
\begin{prop} 
\begin{equation}
(1)\quad D^{\CY}\circ \Rbb=\Rbb\circ D^{\top},\quad C^{\CY}\circ \Rbb=\Rbb\circ C^{\top}.
\end{equation}
(2)\quad $\Xi\cdot e^{-\Theta(u)}\cdot S^{\CY}$ forms a holomorphic flat frame of $H^{\CY}$ w. r. t. the Gauss-Manin connection $\D^{\CY}$.  
\end{prop}

\subsection{Correspondence between real structures}\label{sec4-4.3}\

\

Since in this section we are comparing the real structures defined on each fiber, we can assume without loss of generality that $f_u=f$. 

Let $\{\Gamma_i^{-},i=1,\dots,\mu\}$ be the set of the Lefschetz thimbles constructed in Theorem \ref{Lef-thim}, which forms a basis of $H_{n+1}(\mathbb{C}^{n+2},f^{\leq -\infty},\mathbb{Z})$. Take any $[\gamma]\in H_n(X_f,\mathbb{Z})$. By Lemma \ref{flat-section}, $\{\delta_{2,-t}(\gamma),t\in\mathbb{R}_+\}$ is a flat section and lies in the integral lattice $H_{n+1}(V_{-t},\mathbb{Z})$. Thus there exist integers $\{c_i,i=1,\dots,\mu\}$ such that
\begin{equation}\label{eq:cal-5}
\delta_{2,-t}(\gamma)=\sum_{i=1}^{\mu}c_i[\Gamma^-_i\cap V_{-t}],\quad\forall t\in\mathbb{R}_{+}
\end{equation}

\begin{lm}\label{sec4:lemm-real-stru-2} For any $A\in R_{f}^{(n+2)*}$, the following identity holds:
\begin{equation}\label{sec4:equa-real-1}
\int_\gamma \Rb(A)=-\sum^\mu_{i=1}c_i\int_{\Gamma^-_i}e^f A. 
\end{equation}
\end{lm}

\begin{proof} We can assume that $A$ is a homogeneous element with degree $j$. By Theorem \ref{harmonic-holomorphic-integral} and Remark \ref{coef-discuss}, $\forall i \in \{1,\dots,\mu\}$, we have
\begin{equation}\label{eq:cal-1}
\begin{aligned}\
&\int_{\Gamma^{-}_i}e^{f} A=\int_0^{+\infty}e^{-t}(\int_{\Gamma_i\bigcap V_{-t}} \frac{A dz_1\wedge\dots\wedge dz_{n+2}}{df})dt\\
=& j!\int_{\Gamma_i^{-}\cap V_{-1}}\frac{A dz_1\wedge\dots\wedge dz_{n+2}}{df}\\
=&j!\cdot \frac{1}{2\pi i}\int_{\delta_{1,-1}(\Gamma_i^{-}\cap V_{-1})}\frac{ A dz_1\wedge\dots\wedge dz_{n+2}}{f+1}\\
=&j! Q_i(A,-1),
\end{aligned}
\end{equation}
where the function $Q_{i}(A,t)$ is given by (\ref{sec4:eq-Q-1}). In the proof of \eqref{eq:cal-1}, we have used the fact that $\delta_{1,-1}$ is the Leray coboundary map from $V_{-1}$ to $(\mathbb{C}^{n+2}-V_{-1})$.

By the definition of the map $\Rb$, Theorem \ref{residue-limit}, Lemma \ref{residual-commu}, the definition of the Leray coboundary map $\delta_{2,-t}$, the identity (\ref{eq:cal-5}), and Lemma \ref{residue-integral}, we have
\begin{equation}\label{eq:cal-6}
\begin{aligned}\
&\int_{\gamma}\Rb(A)=(-1)^jj!\int_{\gamma}\Res(\Omega_{A}^{\mathbb{CP}^{n+1}})
=(-1)^jj!\int_{\gamma}\lim_{t\to 0+} r_{4}r_{3,-t}(\omega_{A,-t}) \\
=&(-1)^{j+1} j! \lim_{t\to 0+}\int_{\gamma}r_{2,-t} R_{-t}(\omega_{A,-t})=(-1)^{j+1} j! \lim_{t\to 0+} \int_{\delta_{2,-t}(\gamma)} R_{-t}(\omega_{A,-t})\\
=&-j! \sum_{i=1}^{\mu}c_i Q_{i}(A,-1)=-\sum_{i=1}^{\mu}c_i\int_{\Gamma_i^{-}}e^{f}A.
\end{aligned}
\end{equation}
This finishes the proof of the lemma. 

\end{proof}

\begin{thm}\label{real-structure-correspondence} The real structure $\kappa^{\CY}$ is the usual complex conjugate of the complex vector space $H^{\CY}$, and there is
\begin{align}
\quad \kappa^{\CY}\circ \Rbb=\Rbb\circ \kappa^{\LG}
\end{align}
\end{thm}

\begin{proof} We use the notations in the proof of Lemma \ref{sec4:lemm-real-stru-2}. Let $S^{\CY}_j(u)=\Rb(A_j)$ be a holomorphic local basis of $H^{\CY}$. By Lemma \ref{sec4:lemm-real-stru-2}, we have  
\begin{align*}
&\int_\gamma \kappa^{\CY}\cdot S^{\CY}_j(u)=(\K^{\top})_{\jb}^k \int_\gamma \Rb(A_k)=(\K^{\top})_{\jb}^k \int_{-\sum^\mu_{i=1}c_i \Gamma^-_i} e^f A_k\\
=&\int_{-\sum^\mu_{i=1}c_i \Gamma^-_i} e^f \overline{A_k}=\overline{\int_\gamma \Rb(A_k)}=\int_\gamma \overline{S^{\CY}_j(u)}.
\end{align*}

The commutativity follows easily from the definition of $\kappa^{\CY}$.
\end{proof}

\subsection{Correpondence between pairings}\label{sec4-4.5}

To construct the correspondence between the pairings on the two sides, we need to modify the existed pairings on the CY side and then build the correspondence with the pairing defined in the small $tt^*$ structure on the LG side.

Note that on the LG side, we have the identity
\begin{equation}
\eta^{\top}(S_{\ominus, i}(u),S_{\ominus, j}(u))=\eta^{LG}(\I^{-1}(A_i), \I^{-1}(A_j))(u)=\frac{1}{\mu}\res_{f_u,0}(A_i A_j).
\end{equation}

In particular, for the special pair $[1]$ and $[1^\vee]=[\det(\frac{\pat^2 f}{\pat z_i \pat_{z_j}})]$ in the subring $R_f^{(n+2)*}$ we have
$$
\Res_f([1][1^\vee])=(2\pi i)^n \mu,
$$
where $\mu$ is the Milnor number of $f$.

Now we consider the pairing on the CY side. We start from the following result.

\begin{thm}[{Theorem \ref{fro-iso} or \cite[Theorem 3]{CG}}] For two homogeneous polynomials $A,B\in R_f$ with degrees $\deg A=(n+2)a$ and $\deg B=(n+2)b$ such that $a+b=n$, we have
\begin{equation}\label{sec4:eq-CY-resi-1}
\int_{X_f}(\Res\Omega_A^{\mathbb{CP}^{n+1}})^{n-a,a}\wedge (\Res\Omega_B^{\mathbb{CP}^{n+1}})^{n-b,b}=k_{ab}\Res_{f,0}(A B)
\end{equation}
where
\begin{equation}
k_{ab}=\frac{(-1)^{\frac{a(a+1)+b(b+1)}{2}+b^2+n}}{a!b!}(n+2).
\end{equation}
\end{thm}

Recall the definition of $\Rb$:
$$
\Rb(A)=(-1)^a a! \Res (\Omega_A^{\mathbb{CP}^{n+1}}),
$$
we can rewrite (\ref{sec4:eq-CY-resi-1}) as
\begin{equation*}
\int_{X_f} \Rb(A)\wedge \Rb(B)=i^{n+n(a-b)}(n+2)\Res_{f,0}(AB),
\end{equation*}
or
\begin{equation}
i^{n(b-a-1)}\int_{X_f} \Rb(A)\wedge \Rb(B)=(n+2)\Res_{f,0}(AB).
\end{equation}
In particular, we have
\begin{equation}
i^{n(n-1)}\int_{X_f} \Rb(1)\wedge \Rb(1^\vee)=(n+2)\mu (2\pi i)^{n+2}.
\end{equation}

\begin{df}\label{sec4:defn-CY-pair-1} For any point $u\in M$, we define the pairing on the fiber $H^{\CY}_u$ as
\begin{equation}
\eta^{CY}(S^{\CY}_{A_i}, S^{CY}_{A_j})(u)=\frac{i^{n(\frac{\deg A_j-\deg A_i}{n+2}-1)}\int_{X_{f_u}} \Rb(A_i)\wedge \Rb(A_j)}{i^{n(n-1)}\int_{X_{f_u}} \Rb(1)\wedge \Rb(1^\vee)}=\frac{1}{\mu}\res_{f_u,0}(A_i A_j).
\end{equation}
Then we linearly extend to $H_u^{\CY}$.

the pairing $\eta^{\CY}$ on $H^{\CY}_u$ is defined as the linear extension of the above definition.
\end{df}

\begin{prop}The following conclusions hold.
\begin{enumerate}
\item The pairing $\eta^{\top}$ defined on $\cho\to M$ endows each fiber with a Frobenius algebra structure.
\item The pairing $\eta^{\CY}$ defined on $H^{\CY}\to M$ endows each fiber with a Frobenius algebra structure.
\item The two pairings satisfy the identity:
\begin{equation}
\eta^{\top}=\Rbb^*\eta^{\CY}.
\end{equation}
\end{enumerate}
\end{prop}
\begin{proof} To prove (1) and (2), consider any two homogeneous polynomials $A_i, A_j\in R^{(n+2)*}_{f_u}$ satisfying the degree condition
$$
\frac{\deg A_i+\deg A_j}{n+2}=n.
$$
Let $\pat_1$ represent the deformation direction with respect to $[1]\in R_{f_u}$. We have
$$
\eta^{\top}(S_{\ominus,A_i}, S_{\ominus,A_j})=\frac{1}{\mu}\res_{f_u,0}(A_i A_j)=\eta^{\top}(S_{\ominus,1}, S_{\ominus,A_i A_j}).
$$
We have the same argument for the CY side.

For (3), we have
\begin{align*}
\Rbb^*\eta^{\CY}(S_{\ominus,A_i}, S_{\ominus,A_j})(u)=\eta^{\CY}(S^{\CY}_{A_i}, S^{\CY}_{A_j})(u)=\frac{1}{\mu}\res_{f_u,0}(A_i A_j)=\eta^{\top}(S_{\ominus,A_i}, S_{\ominus,A_j})(u).
\end{align*}
\end{proof}

\subsection{Main theorem}\

\

Let $f:\C^{n+2}\rightarrow \C$ be a nondegenerate homogeneous polynomial of degree $(n+2)$. Let $M\ni 0$ be the marginal deformation space parametrirzed by $R_f^{n+2}\subset R_f$. Then on the CY side we have the new $tt^*$ structure built in Section \ref{sec4-4.4}:
$$
\E^{\CY}=(H^{\CY}\to M, \kappa^{\CY},\eta^{\CY}, D^{\CY}, \bar{D}^{\CY}, C^{\CY}, \bar{C}^{\CY})
$$
with the following data:
\begin{itemize}
\item the fiber $H^{\CY}_u$ at $u\in M$ is the $n$-th primitive cohomology $H^n_{prim}(X_{f_u}, \C)$.
\item $\kappa^{\CY}$ is the real structure on $H^{\CY}$ given by the complex conjugate (see Corollary \ref{real-structure-correspondence}).
\item $\eta^{CY}$ is the pairing given in Definition \ref{sec4:defn-CY-pair-1}.
\item the $tt^*$ connections and the Higgs fields $D^{\CY}, \bar{D}^{\CY}, C^{\CY}, \bar{C}^{\CY}$ are given in Definition \ref{def-connection-higgs}.
\end{itemize}
On the LG side, we have the $tt^*$ substructure given in Definition-Theorem \ref{thm-small-lg}:
$$
\E^{\LG}=({H}^{\LG}\to M, \kappa^{\LG},\eta^{\LG}, D^{\LG}, C^{\LG}, \bar{C}^{\LG}).
$$
Finally, we can summarize our results in Section \ref{Sec4} which, combining Theorem \ref{sec2:theo-harm-coho-corr-1} to give the main theorem in this paper.
\begin{thm}\label{complete-corr}
Let $f:\C^{n+2}\rightarrow \C$ be a nondegenerate homogeneous polynomial of degree $(n+2)$. Then the map $\Rbb\Phi: \E^{\LG}\to \E^{\CY}$ is an isomorphism between $tt^*$ bundles.
\end{thm}

\section{Correspondence for non-Calabi-Yau hypersurfaces}

Let $f:\C^{n+2}\rightarrow \C$ be a nondegenerate homogeneous polynomial of degree $d>1$.$f$ defines a degree $d$ hypersurface $X_f$ in the projective space $\C P^{n+1}$.

The discussion in the above sections can be applied easily to the general (non-Calabi-Yau) hypersurface $X_f$ to obtain partial correspondence between $tt^*$ structure of LG models and $tt^*$ structure of CY model on $X_f$. 

Let $a\in R_f$ be homogeneous element, we define a $\C$-vector space in $R_f$:
%and let $A=adz_1\wedge \cdots \wedge dz_{n+2}$. Then we can define the weight of $a$ and $A$ to be 
%\begin{equation} 
%|a|=\frac{\deg a}{d},\;|A|:=\frac{\deg a+(n+2)}{d}. 
%\end{equation}

\begin{equation}
R_f^{d*}=\{a\in R_f|\frac{\deg a+(n+2)}{d}\in \N\}. 
\end{equation}

\begin{lm}\label{sec5:lm-1} We have the following facts:
\begin{enumerate}
 \item $1^{\vee} \in R_f^{d*}$ for any $d>1$.
\item $1\in R_f^{d*}$ if and only if $d\le (n+2)$ and $d\;|\;(n+2)$. In this case, $R_f^{d*}$ forms a subalgebra of $R_f$ with unit.  
\item If $d>(n+2)$, then $R_f^{d*}$ is not an polynomial algebra with unit. 
\end{enumerate}
\end{lm}

\begin{proof} Note that 
$$
1^\vee=\det(\frac{\pat^2 f}{\pat z_i \pat z_j})=d(d-1)z^{(d-2,\dots,d-2)}\; \text{and}\; |1^\vee|=d-1\in \N.
$$
\end{proof}
Hence $1^\vee\in R_f^{d*}$. 

Let $f_u$ be the deformation of $f$ for $u\in M\subset S_{mr}$. Then $\I^{-1}(R_{f_u}^{d*})$ will induce a subbundle of the big $tt^*$ structure $\hat{\E}^{\LG}$. Simililar to Definition-Theorem \ref{thm-small-lg}, we have 
$$
\E^{\LG}=({H}^{\LG}\to M, \kappa^{\LG},\eta^{\LG}, D^{\LG}, C^{\LG}, \bar{C}^{\LG}),
$$
such that $\I({H}^{\LG}_u)=R_{f_u}^{d*}$. 

On the other hand, the big residue map $\Rbb$ factorizes through the small residue map $\Rb$ such that only the part $R_f^{d*}$ survives and more over $\Rb$ is an isomorphism from $R_{f_u}^{d^*}$ to its image $\Rb(R_{f_u}^{d^*})$ by Lemma \ref{sec4:lemm-real-stru-2}. By a conclusion in \cite{CG}, $\Rb(R_{f_u}^{d*})=H^n (X_{f_u},\C)$. Hence by Definition \ref{def-connection-higgs}, we can construct a (new) $tt^*$ bundle on the moduli space $M$:
$$
\E^{\CY}=(H^{\CY}\to M, \kappa^{\CY},\eta^{\CY}, D^{\CY}, \bar{D}^{\CY}, C^{\CY}, \bar{C}^{\CY}),
$$

\begin{thm}\label{sec5:theo-1} Let $f:\C^{n+2}\rightarrow \C$ be a nondegenerate homogeneous polynomial of degree $d>1$. 
\begin{enumerate}
\item If $d>n+2$, then the (pre)-$tt^*$ bundle structure $(H^{\CY}\to M, \kappa^{\CY},\eta^{\CY}, \D^{CY})$ is just the (pre)-$tt^*$ bundle structure from the classical $tt^*$ bundle $\tilde{\E}^{CY}$. There is no Frobenius algebra structure on the fiber. 

\item If $d<n+2$, then the (pre)-$tt^*$ bundle structure $(H^{\CY}\to M, \kappa^{\CY},\eta^{\CY}, \D^{CY})$ is a (pre)-$tt^*$ subbundle structure from the classical $tt^*$ bundle $\tilde{\E}^{CY}$. There is no Frobenius algebra structure on the fiber if $d\not|\; (n+2)$. 

\item If $d<n+2$ and $d\;|\;(n+2)$, then $\E^{\CY}$ is a $tt^*$ bundle over $M$, and the (pre)-$tt^*$ bundle structure $(H^{\CY}\to M, \kappa^{\CY},\eta^{\CY}, \D^{CY})$ is a (pre)-$tt^*$ subbundle structure from the classical $tt^*$ bundle $\tilde{\E}^{CY}$. There is a Frobenius algebra structure on the fiber of $H^{\CY}$.

\end{enumerate}
\end{thm}

\begin{proof} By (1) of Theorem \ref{prop-residual-map}, to show that the map $\Rb: R_f^{d*}\to H^n_{prim}(X_f,\C)$ is an into or onto map, it suffices to check if any monomial $a$ with $\deg a<d(n+1)$ lies in $R_f^{d*}$. By Lemma \ref{sec5:lm-1}, there is a unique element $1^\vee\in R_f^{d*}$ with the largest degree $\deg 1^\vee=(d-2)(n+2)$. If $d>n+2$, then $\deg 1^\vee+(n+2)>d(n+1)$ and this shows that $\Rb$ is an onto map. If $d<n+2$, then $\deg 1^\vee+(n+2)<d(n+1)$ and this shows that $\Rb$ is an into map. If $d<n+2$ and $d|(n+2)$, then $R_f^{d*}$ forms a subalgebra, and $\E^{\CY}$ is a $tt^*$ subbundle over $M$ since the action of the Higgs fields are closed.  

\end{proof}

\appendix

\section{Gelfand Leray forms and monodromy of Gauss-Manin connection}\label{Appendix-A}

Assume that $f:\mathbb{C}^{n+2}\longrightarrow\mathbb{C}$ is a non-degenerate quasi-homogenous polynomial.\par
Let's consider the singularity theory of $(f,\mathbb{C}^{n+2},0)$. Let $\Delta$ be a small disc near the origin 0, and $\Delta^*=\Delta-\{0\}$. $f$ gives the Milnor fibration $f^{-1}(\Delta^*)\longrightarrow \Delta^*$. Assocaite each $t\in\Delta^*$ with the cohomology group $H^{n+1}(f_t)$, we get the flat vector bundle $H\longrightarrow \Delta^*$ equipped with the Gauss-Manin conncetion. To calculate the monodromy, we should introduce the Brieskorn lattice. For more details about the Brieskorn lattice, see \cite{B} or \cite{Het2}.\par
We know that the non-singular hypersurfaces $f_t$ are stein manifolds, so every cohomology class in $H^{n+1}(f_t)$ can be represented by a holomorphic (n+1)-form on it.\par
There are two ways to get holomorphic (n+1)-forms on the non-singular hypersurfaces $\{f_t\}_{t\in \Delta^*}$. One way is to restrict a holomorphic (n+1)-form on $\mathbb{C}^{n+2}$ to each $f_t$. This gives a subspace of the space of holomrophic sections of $H\longrightarrow\Delta^*$, we shall denote this subspace by $H_f^{\prime}$. \par
Another way is to take the Gelfand-Leray form of a holomorphic $(n+2)$-form. Given a $(n+2)$-form $\omega$ in $\mathbb{C}^{n+2}$, the Gelfand-Leray form of $\omega$ is a holomorphic form in  $H^{n+1}(f_t)$ defined as follow
\begin{displaymath}
\psi(\omega)=\frac{\omega}{df}
\end{displaymath}
As the non-singular hypersurface $f_t$ is given by a regular value $t\in\mathbb{C}^*$, we must have $\omega=df\wedge \theta$ in a neighborhood of $f_t$, and we define $\frac{\omega}{df}$ to be $\theta\mid_{f_t}$. The restriction is independent of the choice of the neighborhood and $\theta$, so $\frac{\omega}{df}$ is a well-defined holomorphic (n+1)-form on $f_t$.\par
Taking the Gelfand-Leray forms defines a subspace of the space of holomorphic forms $H\longrightarrow\Delta^*$. We shall denote this subspace by $H_f^{\prime\prime}$.\par 
We have the following result (see \cite{Seb} and \cite{Het2}):
\begin{thm}
$H_f^{\prime}$ and $H_f^{\prime\prime}$ are both free $\mathcal{O}_{\Delta}$-modules of rank $\mu$. Restrict to the germ at 0$\in\Delta$, we have
\begin{displaymath}
H_{f,0}^{\prime}\cong\Omega^{n+1}_{\mathbb{C}^{n+2},0}/(df\wedge\Omega^n_{\mathbb{C}^{n+2},0}+d\Omega^n_{\mathbb{C}^{n+2},0})
\end{displaymath}
\begin{displaymath}
H_{f,0}^{\prime\prime}\cong\Omega^{n+2}_{\mathbb{C}^{n+2},0}/df\wedge d\Omega^n_{\mathbb{C}^{n+2},0}
\end{displaymath}
And there is a natural embedding $H_f^{\prime}\hookrightarrow H_f^{\prime\prime}$, given by
\begin{displaymath}
[\omega]\longrightarrow [df\wedge\omega]
\end{displaymath}
Consider $H_{f,0}^{\prime}$ as a sub-module of $H_{f,0}^{\prime\prime}$, we have
\begin{displaymath}
H_{f,0}^{\prime\prime}/H_{f,0}^{\prime}\cong\Omega^{n+2}_{\mathbb{C}^{n+2},0}/df\wedge\Omega^{n+1}_{\mathbb{C}^{n+2},0}\cong R_f
\end{displaymath}
\end{thm} 
When $f=F$ is a holomorphic function with the only critical point at the origin 0, by the analytic version of Nullstellensatz, we know that there exists a positive integer $\kappa_F$ such that $F^{\kappa_F}\in (\partial_iF)$. For an arbitrary element $[\omega]\in H_{F,0}^{\prime\prime}$, we have $z^{\kappa_F}[\omega]=[F^{\kappa_F}\omega]\in H_{F,0}^{\prime}$. If $F$ is a quasi-homogeneous polynomial, we can take $\kappa_F=1$. In this case $zH_{F,0}^{\prime\prime}\subseteq  H_{F,0}^{\prime}$, but we know from the above theorem that $H_{F,0}^{\prime\prime}$ is a free $\mathcal{O}_{\Delta,0}$ module of rank $\mu$ and $H_{F,0}^{\prime\prime}/H_{F,0}^{\prime}\cong\Omega^{n+2}_{\mathbb{C}^{n+2},0}/dF\wedge\Omega^{n+1}_{\mathbb{C}^{n+2},0}\cong R_F$ has complex dimension $\mu$, so we must have $zH_{F,0}^{\prime\prime}=H_{F,0}^{\prime}$.\par
We have the following theorem about the action of Gauss-Manin connection on a holomorphic section from $H_F^{\prime}$ \cite{B}.
\begin{thm}
Let $F:\mathbb{C}^{n+2}\longrightarrow\mathbb{C}$ be a holomorphic function with the only critical point at the origin 0. Suppose $\omega$ is a holomorphic (n+1)-form near the origin 0$\in\mathbb{C}^{n+2}$, and $[\omega]$ is the holomorphic section of the flat bundle $H\longrightarrow\Delta^*$ by restricting $\omega$ to non-singular hypersurfaces $\{F_t\}_{t\in\Delta^*}$. Then the action of the Gauss-Manin connection on $[\omega]$ has the form
\begin{displaymath}
\nabla_t[\omega]=[\frac{d\omega}{dF}].
\end{displaymath}
\end{thm}
\begin{proof}
Given an arbitrary $t_0\in\Delta^*$. Let $\gamma (t)$ be a flat section of homology classes $\gamma (t) \in H_n(F_t)$ near $t_0$. We need to calculate the following
\begin{displaymath}
\partial_t\int_{\gamma (t)}\omega
\end{displaymath}
By using the residue formula \cite{Le}, we have 
\begin{displaymath}
\int_{\gamma (t)}\omega=\frac{1}{2\pi i}\int_{\delta\gamma (t)}\frac{dF\wedge\omega}{F-t}
\end{displaymath}
Here $\delta:H_{n+1}(F_t)\longrightarrow H_{n+2}(\mathbb{C}^{n+2}-F_t)$ is the Leray coboundary map, which is defined by taking the boundary of a tubular neighborhood of the homology class in $H_{n+1}(F_t)$. In our case $\gamma (t)$ is a flat section in a small neighborhood of $t_0\in\Delta^*$, $\delta\gamma (t)$ can be taken independent of $t$. Then we have
\begin{displaymath}
\partial_t\int_{\gamma (t)}\omega=\frac{1}{2\pi i}\partial_t\int_{\delta\gamma (t)}\frac{dF\wedge\omega}{F-t}
\end{displaymath}
\begin{displaymath}
=\frac{1}{2\pi i}\int_{\delta\gamma (t)}\frac{dF\wedge\omega}{(F-t)^2}=\frac{1}{2\pi i}(\int_{\delta\gamma (t)}\frac{d\omega}{F-t}-d\frac{\omega}{F-t})
\end{displaymath}
\begin{displaymath}
=\frac{1}{2\pi i}\int_{\delta\gamma (t)}\frac{d\omega}{F-t}=\frac{1}{2\pi i}\int_{\delta\gamma (t)}\frac{dF\wedge \frac{d\omega}{dF}}{F-t}=\int_{\gamma (t)}\frac{d\omega}{dF}
\end{displaymath}
As we can take an arbitrary flat section $\gamma (t)$ near $t_0$, we have proved that 
\begin{displaymath}
\nabla_t[\omega]=[\frac{d\omega}{dF}].
\end{displaymath}
\end{proof}
Given an element $[\omega]\in H_F^{\prime}$ , we can identify it with $[dF\wedge\omega]\in H_F^{\prime\prime}$. Under this embedding, the action of the Gauss-Manin connection has the form
\begin{displaymath}
\nabla_t[dF\wedge \omega]=[d\omega].
\end{displaymath}
In our case, $F=f$ is a (quasi-)homogeneous polynomial, we have $zH_{f,0}^{\prime\prime}=H_{f,0}^{\prime}$. Consider the germ at 0$\in \Delta$, the action of Gauss-Manin connection on $[\beta]\in H_{f,0}^{\prime\prime}$  can be calculate as 
\begin{displaymath}
\nabla_z[\beta]=\nabla_zz^{-1}(z[\beta])=-\frac{[\beta]}{z}+\frac{\nabla_z(z[\beta])}{z}
\end{displaymath}
The last part $\frac{\nabla_z(z[\beta])}{z}$ can be calculated from the formula of Gauss-Manin connection on $H_{f,0}^{\prime}$.\par
Assume that $f(z_1,...,z_{n+2})$ is a quasi-homogeneous polynomial such that 
\begin{displaymath}
f(\lambda^{w_1}z_1,...,\lambda^{w_{n+2}}z_{n+2})=\lambda^{d}f(z_1,...,z_{n+2})
\end{displaymath}
Then we have 
\begin{displaymath}
f(z_1,...,z_{n+2})=\sum_{i=1}^{n+2}\frac{w_i}{d}z_i\partial_{z_i}f(z_1,...,z_{n+2})
\end{displaymath}
We can define a (n+1)-form $\xi$ as
\begin{displaymath}
\xi=\sum_{i=1}^{n+2}(-1)^{i-1}\frac{w_i}{d}z_i dz_1\wedge...\widehat{dz_i}...\wedge dz_{n+2}
\end{displaymath}
Then we have 
\begin{displaymath}
fdz_1\wedge...\wedge dz_{n+2}=df\wedge \xi
\end{displaymath}
With the above discussion, we can prove the following theorem.
\begin{thm}
Assume that $f(z_1,...,z_{n+2})$ is a quasi-homogeneous polynomial with weights $w=(\frac{w_1}{d},...,\frac{w_{n+2}}{d})$. Take a  monomial $\mathbb{C}$-basis $\{z^{\alpha^i}\}_{i=1}^{\mu}$ of $R_f$. Then under the holomorphic basis $\{[z^{\alpha^i}dz_1\wedge...\wedge dz_{n+2}]\}_{i=1}^{\mu}$ of $H_{f,0}^{\prime\prime}$, we have
\begin{displaymath}
z\nabla_z[z^{\alpha^i}dz_1\wedge...\wedge dz_{n+2}]=(<\alpha^i+1,w>-1)[z^{\alpha^i}dz_1\wedge...\wedge dz_{n+2}]
\end{displaymath}
\end{thm}
\begin{proof}
We have 
\begin{displaymath}
z\nabla_z[z^{\alpha^i}dz_1\wedge...\wedge dz_{n+2}]=-[z^{\alpha^i}dz_1\wedge...\wedge dz_{n+2}]+\nabla_zz[z^{\alpha^i}dz_1\wedge...\wedge dz_{n+2}]
\end{displaymath}
\begin{displaymath}
=-[z^{\alpha^i}dz_1\wedge...\wedge dz_{n+2}]+\nabla_z[fz^{\alpha^i}dz_1\wedge...\wedge dz_{n+2}]
\end{displaymath}
\begin{displaymath}
=-[z^{\alpha^i}dz_1\wedge...\wedge dz_{n+2}]+\nabla_z[z^{\alpha^i}df\wedge \xi]=-[z^{\alpha^i}dz_1\wedge...\wedge dz_{n+2}]+[d(z^{\alpha^i}\xi)]
\end{displaymath}
\begin{displaymath}
=-[z^{\alpha^i}dz_1\wedge...\wedge dz_{n+2}]+[\sum_{j=1}^{n+2}\frac{w_j}{d}\partial_j(z^{\alpha^i}\cdot z_j)dz_1\wedge...\wedge dz_{n+2}]
\end{displaymath}
\begin{displaymath}
=(<\alpha^i+1,w>-1)[z^{\alpha^i}dz_1\wedge...\wedge dz_{n+2}]
\end{displaymath}
\end{proof}
When $f$ is a homogeneous polynomial of degree $n+2$, $w=(\frac{1}{n+2},...,\frac{1}{n+2})$, and $<\alpha^i+1,w>-1=<\alpha^i,w>$. So $[z^{\alpha^i}dz_1\wedge...\wedge dz_{n+2}]$ is invariant under the monodromy if and only if deg$z^{\alpha^i}$=$(n+2)k$ for some $k\in\mathbb{Z}_{\geq0}$.\par

\newpage
\bibliographystyle{plain}

\end{document}